\documentclass[10pt]{amsart} 

\usepackage{amsmath}
\usepackage{amssymb}
\usepackage{graphicx}
\usepackage{color}
\usepackage{mathrsfs}
\usepackage{enumitem}
\usepackage{multirow}

\newtheorem{theorem}{Theorem}[section]

\theoremstyle{definition}
\newtheorem{example}[theorem]{Example}
\theoremstyle{remark}

\newtheorem{note}[theorem]{Note}

\newcommand{\R}{\mathbb R}
\newcommand{\C}{\mathcal C}

\newcommand{\A}{\mathscr A}
\newcommand{\B}{\mathscr B}
\newcommand{\D}{\mathscr D}
\newcommand{\E}{\mathscr E}
\newcommand{\F}{\mathscr F}
\newcommand{\G}{\mathscr G}
\newcommand{\HH}{\mathscr H}

\title[Explicit RBF RK methods]{Explicit radial basis function Runge-Kutta methods}
\author{Jiaxi Gu}
\address{Department of Mathematics $\&$ POSTECH MINDS (Mathematical Institute for Data Science), Pohang University of Science and Technology, Pohang 37673, Korea}
\email{jiaxigu@postech.ac.kr}

\author{Xinjuan Chen}
\address{Department of Mathematics, College of Science, Jimei University, Xiamen, Fujian 361021, China}
\email{chenxinjuan@jmu.edu.cn}

\author{Jae-Hun Jung}
\address{Department of Mathematics $\&$ POSTECH MINDS (Mathematical Institute for Data Science), Pohang University of Science and Technology, Pohang 37673, Korea}
\email{jung153@postech.ac.kr}

\makeatletter
\@namedef{subjclassname@2020}{\textup{}2020 Mathematics Subject Classification}
\makeatother

\begin{document}

\maketitle

\begin{abstract}
The aim of this paper is to design the explicit radial basis function (RBF) Runge-Kutta methods for the initial value problem.
We construct the two-, three- and four-stage RBF Runge-Kutta methods based on the Gaussian RBF Euler method with the shape parameter, where the analysis of the local truncation error shows that the $s$-stage RBF Runge-Kutta method could formally achieve order $s\! +\! 1$.
The proof for the convergence of those RBF Runge-Kutta methods follows.
We then plot the stability region of each RBF Runge-Kutta method proposed and compare with the one of the correspondent Runge-Kutta method. 
Numerical experiments are provided to exhibit the improved behavior of the RBF Runge-Kutta methods over the standard ones.
\end{abstract}

\section{Introduction} \label{sec:intro}
The Runge-Kutta methods are widely used for solving the initial value problem of the form
\begin{equation} \label{eq:ivp}
\begin{aligned}
 u' &= f(t,u), \quad a < t \leqslant b \\
 u(a) &= u_0.
\end{aligned}
\end{equation}
The basic idea of the Runge-Kutta methods is to approximate the solution by computing a sequence of intermediate values within a single step in order to achieve high accuracy.
For the explicit $s$-stage Runge-Kutta method with $s \leqslant 4$, we could achieve $s$th-order accuracy.
However, it is impossible to obtain order $s$ with only $s$ stages if $s \geqslant 5$.
For the implicit Runge-Kutta method, there exists some $s$-stage method of order $2s$.
In fact there are numerous papers in the literature on the developments of effective Runge-Kutta methods for specific applications or general classes of problems, including total variation diminishing Runge-Kutta methods \cite{ShuOsherI,GottliebShu} and strong stability preserving Runge-Kutta methods \cite{GottliebShuTadmor,Gottlieb,Bresten} for hyperbolic conservation laws, exponential Runge-Kutta methods \cite{HochbruckSIAM,HochbruckANM,HochbruckActa} for stiff systems of differential equations, and implicit–explicit Runge-Kutta methods \cite{Pareschi,Boscarino,Sebastiano} for time-dependent partial differential equations, among many other.
In \cite{Carpenter,Abarbanel}, the Runge-Kutta methods for the temporal integration of the initial boundary value problem were studied to remove the boundary error and retain the formal order of accuracy in time.

The radial basis function (RBF) interpolation \cite{Buhmann} is one of the most applied approaches in modern approximation theory. 
In \cite{GuJungMQ,GuJungGA}, we derived the multiquadric and Gaussian RBF Euler methods from the RBF interpolation, which outperform the original Euler method in terms of the global error and the order of accuracy. 
The improvement depends on the optimization of the RBF's shape parameter $\epsilon \in \R \cup i \R$ so as to vanish the leading term(s) in the local truncation error. 
Following the idea of generalizing the Euler method to the explicit Runge-Kutta method, we would like to extend the Gaussian RBF Euler method \cite{GuJungGA},
\begin{equation} \label{eq:garbf_euler}
 v_{n+1} = v_n e^{- \epsilon^2_n h^2} + h f(t_n, v_n),
\end{equation}
to the explicit RBF Runge-Kutta method so that the shape parameters are employed to eliminate the leading error term and enhance the formal order by one.
The $s$-stage RBF Runge-Kutta method requires $s-1$ shape parameters in the intermediate stages.
Unlike the identical stability function for all $s$-stage Runge-Kutta methods, the corresponding RBF Runge-Kutta method takes a different form of the stability function due to the choice of shape parameters.
We show that the proposed methods are accurate and efficient in such a way that with a small number of additional operations, at least the overall accuracy is improved via tuning the shape parameters in every step, compared to the standard Runge-Kutta methods. 

The outline of the paper is as follows.  
In Sect. \ref{sec:rk}, we give a brief review of the classical explicit $s$-stage Runge-Kutta methods of order $s$ with $s \leqslant 4$.
Sect. \ref{sec:rbf_rk} defines a class of RBF Runge-Kutta methods up to four stages and provides conditions under which the local truncation error of the $s$-stage RBF Runge-Kutta method is of order $s\! +\! 1$.
Our convergence analysis of the RBF Runge-Kutta methods is performed in Sect. \ref{sec:convg}.
The plot of the stability regions is contained in Sect. \ref{sec:sregion}, where we compare the stability regions of the standard and RBF Runge-Kutta methods with the same number of stages.
In Sect. \ref{sec:nr} we present numerical experiments which show that (i) it is possible for the $s$-stage RBF Runge-Kutta method to attain $(s\! +\! 1)$th-order accuracy, (ii) the order increase/reduction for some RBF Runge-Kutta method arises in practical examples and (iii) the $s$-stage RBF Runge-Kutta method obtains better accuracy than the $s+1$-stage Runge-Kutta method for some problems.
Concluding remarks are given in Sect. \ref{sec:conclusion}.

\section{Explicit Runge-Kutta methods} \label{sec:rk}
The general explicit $s$-stage Runge-Kutta methods are given by 
$$ 
   v_{n+1} = v_n + h \sum_{i=1}^s b_i k_i,
$$
where $h$ is the step size and
\begin{align*}
 k_1 &= f(t_n, v_n), \\
 k_2 &= f(t_n + c_2 h, v_n + h (a_{21} k_1)), \\
 k_3 &= f(t_n + c_3 h, v_n + h (a_{31} k_1 + a_{32} k_2)), \\
 \vdots & \\
 k_s &= f \left( t_n + c_s h, v_n + h \sum_{j=1}^{s-1} a_{ij} k_j \right).
\end{align*}
Here $a_{ij} \left( 1 \leqslant j < i \leqslant s \right), \, b_i \left( i = 1, \cdots, s \right)$ and $c_i \left( i = 2, \cdots, s \right)$ are to be determined. 
Each Runge-Kutta method could be represented uniquely by its Butcher tableau \cite{Butcher}
$$ 
\begin{array}{c|cccccr}
   0   &        &        &        &            &     & \\
   c_2 & a_{21} &        &        &            &     & \\
   c_3 & a_{31} & a_{32} &        &            &     & \\
\vdots & \vdots & \vdots & \ddots &            &     & \\
   c_s & a_{s1} & a_{s2} & \cdots & a_{s(s-1)} &     & \: \underset{.}{} \\
\cline{1-6}
       & b_1    & b_2    & \cdots & b_{s-1}    & b_s & 
\end{array}
$$
A Taylor series expansion shows that the Runge-Kutta method is consistent if and only if
\begin{equation}  \label{eq:rk_consistent}
 \sum_{j=1}^{s} b_j = 1.
\end{equation}
Throughout this paper, we assume that the Runge-Kutta method satisfies 
\begin{equation}  \label{eq:rk_cond}
 \sum_{j=1}^{i-1} a_{ij} = c_i, \, i = 1, \cdots, s,
\end{equation}
in order to achieve high order.

\subsection{Two-stage second-order methods}
The two-stage Runge-Kutta method takes the form 
\begin{equation} \label{eq:rk2}
\begin{aligned}
 v_{n+1} &= v_n + h (b_1 k_1 + b_2 k_2), \\
 k_1 &= f(t_n, v_n), \\
 k_2 &= f(t_n + c_2 h, v_n + h a_{21} k_1),
\end{aligned}
\end{equation}
and the corresponding Butcher tableau is
$$
\begin{array}{c|ccr}
 0   &        &     & \\
 c_2 & a_{21} &     & \: \underset{.}{} \\
\cline{1-3}
     & b_1    & b_2 &
\end{array}
$$
By the Taylor series expansion, the Runge-Kutta method \eqref{eq:rk2} is second-order accurate if we have
\begin{equation} \label{eq:rk2_cond}
 b_2 c_2 = \frac{1}{2},
\end{equation}
along with the conditions \eqref{eq:rk_consistent} and \eqref{eq:rk_cond}.
This gives a one-parameter family of two-stage Runge-Kutta methods
\begin{equation} \label{tab:rk22}
\begin{array}{c|cc}
 0   &                  & \\
 c_2 & c_2              & \\
\hline 
     & 1-\frac{1}{2c_2} & \frac{1}{2c_2}
\end{array}
\end{equation}
for $c_2 \ne 0$.
Note that the two-stage second-order Runge-Kutta method \eqref{tab:rk22} is equivalent to
\begin{align*}
 v^{(1)} &= v_n + c_2 h f(t_n, v_n), \\
 v_{n+1} &= \left[ 1 - \frac{1}{c_2} \left( 1 - \frac{1}{2c_2} \right) \right] v_n + \frac{1}{c_2} \left( 1 - \frac{1}{2c_2} \right) v^{(1)} + \frac{1}{2c_2} h f \left( t_n + c_2 h, v^{(1)} \right).
\end{align*}

\subsection{Three-stage third-order methods}
The form of the three-stage Runge-Kutta method is 
\begin{equation} \label{eq:rk3}
\begin{aligned}
 v_{n+1} &= v_n + h (b_1 k_1 + b_2 k_2 + b_3 k_3), \\
 k_1 &= f(t_n, v_n), \\
 k_2 &= f(t_n + c_2 h, v_n + h a_{21} k_1), \\
 k_3 &= f(t_n + c_3 h, v_n + h (a_{31} k_1 + a_{32} k_2)).
\end{aligned}
\end{equation}
Expressed in a Butcher tableau, the method becomes
$$
\begin{array}{c|cccr}
 0   &        &        &     & \\
 c_2 & a_{21} &        &     & \\
 c_3 & a_{31} & a_{32} &     & \: \underset{.}{} \\
\cline{1-4}
     & b_1    & b_2    & b_3 & 
\end{array} 
$$
The three-stage Runge-Kutta method \eqref{eq:rk3} can also take the alternative form
\begin{align*}
 v^{(1)} &= v_n + c_2 h f(t_n, v_n), \\
 v^{(2)} &= \left( 1-\frac{a_{31}}{c_2} \right) v_n + \frac{a_{31}}{c_2} v^{(1)} + a_{32} h f \left( t_n+c_2 h, v^{(1)} \right), \\
 v_{n+1} &= \left[ 1 - \frac{1}{c_2} \left( b_1 - \frac{a_{31}}{a_{32}} b_2 \right) - \frac{b_2}{a_{32}} \right] v_n + 
            \frac{1}{c_2} \left( b_1 - \frac{a_{31}}{a_{32}} b_2 \right) v^{(1)} + \frac{b_2}{a_{32}} v^{(2)} + 
            b_3 h f \left( t_n + c_3 h, v^{(2)} \right).
\end{align*}
The conditions for order $3$ are
\begin{equation} \label{eq:rk3_cond}
\begin{aligned} 
 b_2 c_2 + b_3 c_3 &= \frac{1}{2}, \\
 b_2 c^2_2 + b_3 c^2_3 &= \frac{1}{3}, \\
 a_{32} b_3 c_2 &= \frac{1}{6},
\end{aligned}
\end{equation}
in combination with the conditions \eqref{eq:rk_consistent} and \eqref{eq:rk_cond}.
It is shown in \cite{Butcher} that there are three families of three-stage third-order Runge-Kutta methods ($a_{32} \ne 0$):
\begin{enumerate}[label=\Roman*.]
\item $c_2 \ne 0, \, c_2 \ne \frac{2}{3}, \, c_3 \ne 0, \, c_2 \ne c_3$ \\
$$
\begin{array}{c|cccr}
 0   &     &        & & \\
 c_2 & c_2 &        & & \\
 c_3 & \frac{c_3 \left( 3c_3-3c_2^2-c_3 \right)}{c_2 \left( 2-3c_2 \right)} & \frac{c_3 \left( c_3-c_2 \right)}{c_2 \left( 2-3c_2 \right)} & & \: \underset{;}{} \\
\cline{1-4}
     & \frac{\left( -3c_3+6c_2c_3+2-3c_3 \right)}{6c_2 c_3} & \frac{\left( 3c_3-2 \right)}{6c_2 \left( c_3-c_2 \right)} & \frac{\left( 2-3c_2 \right)}{6c_3 \left( c_3-c_2 \right)} &
\end{array}
$$
\item $b_3 \ne 0, \, c_2 = c_3 = \frac{2}{3}$ \\
$$
\begin{array}{c|cccr}
 0           &                              &                 &     & \\
 \frac{2}{3} & \frac{2}{3}                  &                 &     & \\
 \frac{2}{3} & \frac{2}{3} - \frac{1}{4b_3} & \frac{1}{4b_3}  &     & \: \underset{;}{} \\
\cline{1-4}
             & \frac{1}{4}                  & \frac{3}{4}-b_3 & b_3 & 
\end{array}
$$
\item $b_3 \ne 0, \, c_2 = \frac{2}{3}, \, c_3 = 0$ \\
$$
\begin{array}{c|cccr}
 0           &                  &                &     & \\
 \frac{2}{3} & \frac{2}{3}      &                &     & \\
 0           & - \frac{1}{4b_3} & \frac{1}{4b_3} &     & \: \underset{.}{} \\
\cline{1-4}
             & \frac{1}{4}-b_3  & \frac{3}{4}    & b_3 &
\end{array}
$$
\end{enumerate}

\subsection{Four-stage fourth-order methods}
The four-stage Runge-Kutta method is of form 
\begin{equation} \label{eq:rk4}
\begin{aligned}
 v_{n+1} &= v_n + h (b_1 k_1 + b_2 k_2 + b_3 k_3 + b_4 k_4), \\
 k_1 &= f(t_n, v_n), \\
 k_2 &= f(t_n + c_2 h, v_n + h a_{21} k_1), \\
 k_3 &= f(t_n + c_3 h, v_n + h (a_{31} k_1 + a_{32} k_2)), \\
 k_4 &= f(t_n + c_4 h, v_n + h (a_{41} k_1 + a_{42} k_2 + a_{43} k_3)), \\
\end{aligned}
\end{equation}
or
\begin{align*}
 v^{(1)} = {} & v_n + a_{21} h f(t_n, v_n), \\
 v^{(2)} = {} & \left( 1 - \frac{a_{31}}{a_{21}} \right) v_n + \frac{a_{31}}{a_{21}} v^{(1)} + 
                a_{32} h f \left( t_n + c_2 h, v^{(1)} \right), \\
 v^{(3)} = {} & \left[ 1 - \frac{a_{41}}{a_{21}} - \frac{a_{42}}{a_{32}} \left( 1-\frac{a_{31}}{a_{21}} \right) \right] v_n + 
                \left( \frac{a_{41}}{a_{21}} -  \frac{a_{42}a_{31}}{a_{32}a_{21}} \right) v^{(1)} + \frac{a_{42}}{a_{32}} v^{(2)}+\\
           {} & a_{43} h f \left( t_n + c_3 h, v^{(2)} \right), \\
 v_{n+1} = {} & \left\{ 1 - \frac{b_3}{a_{43}} - \frac{1}{a_{32}} \left( b_2 - b_3 \frac{a_{42}}{a_{43}} \right) - \frac{1}{a_{21}} \left[ b_1 - b_3 \frac{a_{41}}{a_{43}} - \frac{a_{31}}{a_{32}} \left( b_2-b_3 \frac{a_{42}}{a_{43}} \right) \right] \right\} v_n +\\
           {} & \frac{1}{a_{21}}\left[ b_1 - b_3\frac{a_{41}}{a_{43}} - \frac{a_{31}}{a_{32}} \left(b_2 - b_3 \frac{a_{42}}{a_{43}} \right) \right] v^{(1)} + \\
           {} & \frac{1}{a_{32}} \left( b_2 - b_3 \frac{a_{42}}{a_{43}} \right) v^{(2)} + \frac{b_3}{a_{43}} v^{(3)} + 
                b_4 h f \left( t_n + c_4 h, v^{(3)} \right).
\end{align*}
The Butcher tableau of the method is
$$
\begin{array}{c|ccccr}
 0   &        &        &        &     & \\
 c_2 & a_{21} &        &        &     & \\
 c_3 & a_{31} & a_{32} &        &     & \\
 c_4 & a_{41} & a_{42} & a_{43} &     & \: \underset{.}{} \\
\cline{1-5} 
     & b_1    & b_2    & b_3    & b_4 &
\end{array}
$$
Fourth-order accuracy requires the conditions
\begin{equation} \label{eq:rk4_cond}
\begin{aligned}
 b_2 c_2 + b_3 c_3 + b_4 c_4 &= \frac{1}{2}, \\
 b_2 c^2_2 + b_3 c^2_3 + b_4 c^2_4 &= \frac{1}{3}, \\
 b_2 c^3_2 + b_3 c^3_3 + b_4 c^3_4 &= \frac{1}{4}, \\
 a_{32} b_3 c_2 + a_{42} b_4 c_2 + a_{43} b_4 c_3 &= \frac{1}{6}, \\
 a_{32} b_3 c_2 c_3 + a_{42} b_4 c_2 c_4 + a_{43} b_4 c_3 c_4 &= \frac{1}{8}, \\
 a_{32} b_3 c^2_2 + a_{42} b_4 c^2_2 + a_{43} b_4 c^2_3 &= \frac{1}{12}, \\
 a_{32} a_{43} b_4 c_2 &= \frac{1}{24},
\end{aligned}
\end{equation} 
together with conditions \eqref{eq:rk_consistent} and \eqref{eq:rk_cond}.
It is easy to show that $c_4 = 1$.
There is an order barrier for the explicit Runge-Kutta methods, i.e., it is not feasible to use $s$ stages for order $s$ if $s \geqslant 5$, as explained in \cite{Butcher}.
For example, at least six stages are required for order $5$, and seven stages for order $6$.

\section{Explicit RBF Runge-Kutta methods} \label{sec:rbf_rk}
Inspired by the Gaussian RBF Euler method \eqref{eq:garbf_euler}, we now construct the intermediate stage by introducing the shape parameter $\epsilon \in \R \cup i \R$, yielding the explicit RBF Runge-Kutta methods.

\subsection{Two-stage third-order methods} \label{sec:rbf_rk2}
The two-stage RBF Runge-Kutta method is devised to take the form
\begin{equation} \label{eq:rbf_rk2}
\begin{aligned}
 v_{n+1} &= v_n + h (b_1 k_1 + b_2 k_2), \\
     k_1 &= f(t_n, v_n), \\
     k_2 &= f(t_n + c_2 h, v_n e^{- \epsilon^2_n (a_{21} h)^2} + h a_{21} k_1),
\end{aligned}
\end{equation}
which can be expressed by the augmented Butcher tableau
\begin{equation} \label{tab:rbf_rk2}
\begin{array}{c|cc|c}
 0   &        &     & \\
 c_2 & a_{21} &     & \epsilon^2_n \\
\hline 
     & b_1    & b_2 & 
\end{array}
\end{equation}
with the information of $\epsilon^2_n$ appended to the right besides the values of $a_{ij}, \, b_i$ and $c_i$.
This method reduces to the two-stage Runge-Kutta method \eqref{eq:rk2} if $\epsilon^2_n = 0$.  
Similarly, we can rephrase this two-stage method \eqref{eq:rbf_rk2} to
\begin{align*}
 v^{(1)} &= v_n e^{- \epsilon^2_n (a_{21} h)^2} + a_{21} h f(t_n, v_n), \\
 v_{n+1} &= \left[ 1 - \frac{b_1}{a_{21}} e^{- \epsilon^2_n (a_{21} h)^2} \right] v_n + 
            \frac{b_1}{a_{21}} v^{(1)} + b_2 h f(t_n + c_2 h, v^{(1)}).
\end{align*}

By Taylor series expansion, the local truncation error, defined as the error occurring from a single step normalized by $h$ while taking $v_n = u_n$ as the exact solution, is
\begin{align} 
 \tau_n = {} & \frac{u_{n+1} - u_n}{h} - (b_1 k_1 + b_2 k_2) \label{eq:rbf_rk2_lte} \\
        = {} & \left( 1 - b_1 - b_2 \right) f + 
               \left[ \left( \frac{1}{2} - b_2 c_2 \right) f_t + 
                      \left( \frac{1}{2} - a_{21} b_2 \right) f_u f \right] h + \nonumber \\
          {} & \bigg[ \left( \frac{1}{6} - \frac{1}{2} b_2 c^2_2 \right) f_{tt} + 
                      \left( \frac{1}{3} - a_{21} b_2 c_2 \right) f_{tu} f + 
                      \left( \frac{1}{6} - \frac{1}{2} a^2_{21} b_2 \right) f_{uu} f^2 + \nonumber \\
          {} & \left( \frac{f_t + f_u f}{6} + \epsilon^2_n a^2_{21} b_2 u_n \right) f_u \bigg] h^2 + O(h^3). \nonumber
\end{align}
When the conditions \eqref{eq:rk_consistent}, \eqref{eq:rk_cond}, and \eqref{eq:rk2_cond} are satisfied, the two-stage RBF Runge-Kutta method maintains second-order accuracy as the leading error term is of $O(h^2)$.
If we further require 
$$
   b_2 c^2_2 = \frac{1}{3},
$$
then $a_{21} = c_2 = \frac{2}{3}, \, b_1 = \frac{1}{4}, \, b_2 = \frac{3}{4}$, and the first three terms in $h^2$ vanish.
Setting the last term in $h^2$ to $0$ gives
\begin{equation} \label{eq:rbf_rk2_eps}
 \epsilon^2_n = - \frac{f_t + f_u f}{2 u_n} = - \frac{u''_n}{2 u_n}.
\end{equation}
So this RBF Runge-Kutta method is third-order accurate.
\begin{note}
 When $\epsilon^2_n = 0$, it becomes the Ralston's method \cite{Ralston}, which is a second-order method with a minimum local error bound.
\end{note} 

\subsection{Three-stage fourth-order methods} \label{sec:rbf_rk3}
We design the three-stage RBF Runge-Kutta method of form 
\begin{equation} \label{eq:rbf_rk3}
\begin{aligned}
 v_{n+1} &= v_n + h (b_1 k_1 + b_2 k_2 + b_3 k_3), \\
     k_1 &= f(t_n, v_n), \\
     k_2 &= f(t_n + c_2 h, v_n e^{- \epsilon^2_{n2} (c_2 h)^2} + h a_{21} k_1), \\
     k_3 &= f(t_n + c_3 h, v_n e^{- \epsilon^2_{n3} (c_3 h)^2} + h (a_{31} k_1 + a_{32} k_2)),
\end{aligned}
\end{equation}
and its augmented Butcher tableau is 
$$
\begin{array}{c|ccc|cr}
 0   &        &        &     &                 & \\
 c_2 & a_{21} &        &     & \epsilon^2_{n2} & \\
 c_3 & a_{31} & a_{32} &     & \epsilon^2_{n3} & \: \underset{.}{} \\
\cline{1-5} 
     & b_1    & b_2    & b_3 &                 &
\end{array}
$$
The alternative form of this RBF Runge-Kutta methhod is given by
\begin{align*}
 v^{(1)} = {} & v_n e^{- \epsilon^2_{n2} (c_2 h)^2} + a_{21} h f(t_n, v_n), \\
 v^{(2)} = {} & \left[ e^{- \epsilon^2_{n3} (c_3 h)^2} - \frac{a_{31}}{a_{21}} e^{- \epsilon^2_{n2} (c_2 h)^2} \right] v_n +
                \frac{a_{31}}{a_{21}} v^{(1)} + a_{32} h f( t_n + c_2 h, v^{(1)}), \\
 v_{n+1} = {} & \left[ 1 - \frac{1}{a_{21}} \left( b_1 - \frac{a_{31}}{a_{32}} b_2 \right) e^{- \epsilon^2_{n2} (c_2 h)^2} - \frac{b_2}{a_{32}} e^{- \epsilon^2_{n3} (c_3 h)^2} \right] v_n + \\
           {} & \frac{1}{a_{21}} \left( b_1 - \frac{a_{31}}{a_{32}} b_2 \right) v^{(1)} + \frac{b_2}{a_{32}} v^{(2)} + b_3 h f(t_n + c_3 h, v^{(2)}).
\end{align*}

With Taylor series expansion and the conditions $a_{21} = c_2$ and $a_{31} + a_{32} = c_3$ in \eqref{eq:rk_cond}, the local truncation error could be simplified as
\begin{align} 
 \tau_n ={} & \frac{u_{n+1} - u_n}{h} - (b_1 k_1 + b_2 k_2 + b_3 k_3) \label{eq:rbf_rk3_lte} \\
 \begin{split} \label{eq:rbf_rk3_lte_expd}
         = {} & \left( 1 - b_1 - b_2 - b_3 \right) f + 
                \left( \A^3_1 f_t + \A^3_1 f_u f \right) h + 
                \Big[ \A^3_2 f_{tt} + 2 \A^3_2 f_{tu} f + \\ 
           {} & \B^3_2 f_t f_u + \B^3_2 f^2_u f + \A^3_2 f_{uu} f^2 + 
                \left( b_2 c^2_2 \epsilon^2_{n2} + b_3 c^2_3 \epsilon^2_{n3} \right) f_u u_n \Big] h^2 + \\
           {} & \bigg[ \A^3_3 f_{ttt} + 3 \A^3_3 f_{ttu} f + \B^3_3 f_t f_{tu} + 3 \A^3_3 f_{tuu} f^2 + \D^3_3 f_{tt} f_u + \\
           {} & \left( \B^3_3 + 2 \D^3_3 \right) f_{tu} f_u f + \B^3_3 f_t f_{uu} f + 
                \left( \B^3_3 + \D^3_3 \right) f_u f_{uu} f^2 + \A^3_3 f_{uuu} f^3 + \\
           {} & \left( b_2 c^3_2 \epsilon^2_{n2} + b_3 c^3_3 \epsilon^2_{n3} \right) \left( f_{tu} u_n + f_{uu} f u_n \right) +
                \epsilon^2_{n2} a_{32} b_3 c^2_2 f^2_u u_n + \frac{1}{24} f_u^2 u''_n \bigg] h^3 + O(h^4),
 \end{split}
\end{align}
where $\A^3_1, \, \A^3_2, \, \B^3_2, \, \A^3_3, \, \B^3_3, \, \D^3_3$ are given in Appendix.

Applying the conditions \eqref{eq:rk_consistent}, \eqref{eq:rk3_cond} and
$$
   b_2 c^2_2 \epsilon^2_{n2} + b_3 c^2_3 \epsilon^2_{n3} = 0,
$$
to the local truncation error \eqref{eq:rbf_rk3_lte_expd} makes the constant, $h$ and $h^2$ terms vanish.
Thus $\tau_n = O(h^3)$ and this three-stage RBF Runge-Kutta method retains order $3$.
Furthermore, in order to obtain fourth-order accuracy, two additional conditions are required in the following $4$ different cases:
\begin{enumerate}[label=\Roman*.]
\item $\A^3_3 = 0, \, \D^3_3 = 0$ \\ 
      The coefficients of $f_{ttt}, \, f_{ttu} f, \, f_{tuu} f^2, \, f_{tt} f_u$ and $f_{uuu} f^3$ in $h^3$ become zero.
      Then $a_{21} = c_2 = \frac{1}{2}, \, a_{31} = -1, \, a_{32} = 2, \, c_3 = 1, \, b_1 = \frac{1}{6}, \, b_2 = \frac{2}{3}$ and $b_3 = \frac{1}{6}$.       
      Equating the rest of terms in $h^3$ to zero, we have
$$
   \epsilon^2_{n2} = - \frac{u''_n}{2 u_n}, \quad \epsilon^2_{n3} = - \epsilon^2_{n2}.
$$
\item $\A^3_3 = 0, \, \B^3_3 + 2 \D^3_3 = 0$ \\ 
      The $f_{ttt}, \, f_{ttu} f, \, f_{tuu} f^2, \, f_{tu} f_u f$ and $f_{uuu} f^3$ terms in $h^3$ are annihilated.
 \begin{enumerate}
 \item $a_{21} = c_2 = \frac{15-\sqrt{33}}{24}, \, a_{31} = -\frac{147+29\sqrt{33}}{768}, \, a_{32} = \frac{627+61\sqrt{33}}{768}, \, c_3 = \frac{15+\sqrt{33}}{24}, \, b_1 = \frac{1}{8}, \, b_2 = \frac{77+3\sqrt{33}}{176}$ and $b_3 = \frac{77-3\sqrt{33}}{176}$.       
      If we set the leading error term of $O(h^3)$ to zero, then
\begin{align*}
 \epsilon^2_{n2} &= \frac{-2(3-\sqrt{33}) (f_{tu}+f_{uu}f) f_t + (3-\sqrt{33}) (f_{tt}-f_{uu}f^2) f_u - 12 f_u^2 u''_n}{2 [ 2(3-\sqrt{33}) (f_{tu}+f_{uu}f) + (15-\sqrt{33})f_u^2 ] u_n}, \\
 \epsilon^2_{n3} &= - \frac{7-\sqrt{33}}{4} \epsilon^2_{n2}.
\end{align*}  
 \item $a_{21} = c_2 = \frac{15+\sqrt{33}}{24}, \, a_{31} = \frac{-147+29\sqrt{33}}{768}, \, a_{32} = \frac{627-61\sqrt{33}}{768}, \, c_3 = \frac{15-\sqrt{33}}{24}, \, b_1 = \frac{1}{8}, \, b_2 = \frac{77-3\sqrt{33}}{176}$ and $b_3 = \frac{77+3\sqrt{33}}{176}$.       
      If we vanish the leading error term in $O(h^3)$, then
\begin{align*}
 \epsilon^2_{n2} &= \frac{-2(3+\sqrt{33}) (f_{tu}+f_{uu}f) f_t + (3+\sqrt{33}) (f_{tt}-f_{uu}f^2) f_u - 12 f_u^2 u''_n}{2 [ 2(3+\sqrt{33}) (f_{tu}+f_{uu}f) + (15+\sqrt{33})f_u^2 ] u_n}, \\
 \epsilon^2_{n3} &= - \frac{7+\sqrt{33}}{4} \epsilon^2_{n2}.
\end{align*}  
 \end{enumerate}
\item $\A^3_3 = 0, \, \B^3_3 + \D^3_3 = 0$ \\ 
      The coefficients of $f_{ttt}, \, f_{ttu} f, \, f_{tuu} f^2, \, f_u f_{uu} f^2$ and $f_{uuu} f^3$ in $h^3$ are zero.
 \begin{enumerate}
 \item $a_{21} = c_2 = \frac{1}{3}, \, a_{31} = -\frac{5}{12}, \, a_{32} = \frac{5}{4}, \, c_3 = \frac{5}{6}, \, b_1 = \frac{1}{10}, \, b_2 = \frac{1}{2}$ and $b_3 = \frac{2}{5}$.       
       Setting the rest of terms in $h^3$ to zero results in
\begin{align*}
 \epsilon^2_{n2} &= \frac{(f_{tu}+f_{uu}f) f_t - (f_{tt}+f_{tu}f) f_u - 3 f_u^2 u''_n}{2 (2f_u^2-f_{tu}-f_{uu}f) u_n}, \\
 \epsilon^2_{n3} &= - \frac{1}{5} \epsilon^2_{n2}.
\end{align*}  
 \item $a_{21} = c_2 = 1, \, a_{31} = \frac{1}{4}, \, a_{32} = \frac{1}{4}, \, c_3 = \frac{1}{2}, \, b_1 = \frac{1}{6}, \, b_2 = \frac{1}{6}$ and $b_3 = \frac{2}{3}$.       
       Annihilating the leading error term in $O(h^3)$, we get
\begin{align*}
 \epsilon^2_{n2} &= \frac{-(f_{tu}+f_{uu}f) f_t + (f_{tt}+f_{tu}f) f_u - f_u^2 u''_n}{2 (2f_u^2+f_{tu}+f_{uu}f) u_n}, \\
 \epsilon^2_{n3} &= - \epsilon^2_{n2}.
\end{align*}  
 \end{enumerate}
\item $\B^3_3 = 0, \, \D^3_3 = 0$ \\ 
      The $f_t f_{tu}, \, f_{tt} f_u, \, f_{tu} f_u f, \, f_t f_{uu} f$ and $f_u f_{uu} f^2$ terms in $h^3$ vanish.
      Then $a_{21} = c_2 = \frac{1}{2}, \, a_{31} = 0, \, a_{32} = c_3 = \frac{3}{4}, \, b_1 = \frac{2}{9}, \, b_2 = \frac{1}{3}$ and $b_3 = \frac{4}{9}$.       
      If we have the leading error term of $O(h^3)$ equal zero, then
\begin{align*}
 \epsilon^2_{n2} &= - \frac{f_{ttt} + f_{uuu}f^3 + 3 (f_{ttu} + f_{tuu}f)f + 12 f_u^2 u''_n}{6 (4f_u^2-f_{tu}-f_{uu}f) u_n}, \\
 \epsilon^2_{n3} &= - \frac{1}{3} \epsilon^2_{n2}.
\end{align*}  
\end{enumerate}
\begin{note} 
 If we set $\epsilon^2_{n2} = \epsilon^2_{n3} = 0$, the above RBF Runge-Kutta methods return to three-stage third-order Runge-Kutta methods.
\end{note} 

\subsection{Four-stage fifth-order methods} \label{sec:rbf_rk4}
We have seen that introducing a shape parameter in each stage could increase the order of accuracy by one analytically.
This technique could be applied to construct the four-stage RBF Runge-Kutta method,
\begin{equation} \label{eq:rbf_rk4}
\begin{aligned}
 v_{n+1} &= v_n + h (b_1 k_1 + b_2 k_2 + b_3 k_3 + b_4 k_4), \\
     k_1 &= f(t_n, v_n), \\
     k_2 &= f(t_n + c_2 h, v_n e^{- \epsilon^2_{n2} (c_2 h)^2} + h a_{21} k_1), \\
     k_3 &= f(t_n + c_3 h, v_n e^{- \epsilon^2_{n3} (c_3 h)^2} + h (a_{31} k_1 + a_{32} k_2)), \\
     k_4 &= f(t_n + c_4 h, v_n e^{- \epsilon^2_{n4} (c_4 h)^2} + h (a_{41} k_1 + a_{42} k_2 + a_{43} k_3)) ,
\end{aligned}
\end{equation}
which corresponds to the augmented Butcher tableau
$$
\begin{array}{c|cccc|cr}
 0   &        &        &        &     &                 & \\
 c_2 & a_{21} &        &        &     & \epsilon^2_{n2} & \\
 c_3 & a_{31} & a_{32} &        &     & \epsilon^2_{n3} & \\
 c_4 & a_{41} & a_{42} & a_{43} &     & \epsilon^2_{n3} & \: \underset{.}{} \\
\cline{1-6} 
     & b_1    & b_2    & b_3    & b_4 &                 &
\end{array}
$$
The RBF Runge-Kutta method can be rewritten as 
\begin{align*}
 v^{(1)} = {} & v_n e^{-\epsilon^2_{n2}(c_2 h)^2} + a_{21} h f(t_n, v_n), \\
 v^{(2)} = {} & \left( e^{-\epsilon^2_{n3}(c_3 h)^2} - \frac{a_{31}}{a_{21}}e^{-\epsilon^2_{n2}(c_2 h)^2} \right) v_n +
                \frac{a_{31}}{a_{21}} v^{(1)} + a_{32} h f(t_n + c_2 h, v^{(1)}), \\
 v^{(3)} = {} & \left[ e^{-\epsilon^2_{n4}(c_4 h)^2} - \frac{a_{42}}{a_{32}} e^{-\epsilon^2_{n3}(c_3 h)^2} - \frac{1}{a_{21}} \left( a_{41}-\frac{a_{42}}{a_{32}} a_{31} \right) e^{-\epsilon^2_{n2}(c_2 h)^2} \right] v_n + \\
           {} & \frac{1}{a_{21}} \left( a_{41} - \frac{a_{42}}{a_{32}}a_{31} \right) v^{(1)} + \frac{a_{42}}{a_{32}} v^{(2)} +
                a_{43} h f(t_n + c_3 h, v^{(2)}), \\
 v_{n+1} = {} & \left\{ 1 - \frac{b_3}{a_{43}} e^{-\epsilon^2_{n4}(c_4 h)^2} - \frac{1}{a_{32}} \left( b_2 - b_3 \frac{a_{42}}{a_{43}} \right) e^{-\epsilon^2_{n3}(c_3 h)^2} - \right. \\
           {} & \left. \frac{1}{a_{21}} \left[ b_1 - b_3 \frac{a_{41}}{a_{43}} - \frac{a_{31}}{a_{32}} \left( b_2-b_3 \frac{a_{42}}{a_{43}} \right) \right] e^{-\epsilon^2_{n2}(c_2 h)^2} \right\} v_n + \\
           {} & \frac{1}{a_{21}} \left[ b_1 - b_3 \frac{a_{41}}{a_{43}} - \frac{a_{31}}{a_{32}} \left( b_2 - b_3 \frac{a_{42}}{a_{43}} \right) \right] v^{(1)} + \frac{1}{a_{32}} \left( b_2 - b_3 \frac{a_{42}}{a_{43}} \right) v^{(2)} + \\
           {} & \frac{b_3}{a_{43}} v^{(3)} + b_4 h f(t_n + c_4 h, v^{(3)}).
\end{align*}

Taylor series expansion, with the conditions $a_{21} = c_2, \, a_{31} + a_{32} = c_3, \, a_{41} + a_{42} + a_{43} = c_4$ in \eqref{eq:rk_cond}, yields 
\begin{align}
 \tau_n = {} & \frac{u_{n+1} - u_n}{h} - (b_1 k_1 + b_2 k_2 + b_3 k_3 + b_4 k_4) \label{eq:rbf_rk4_lte} \\
 \begin{split} \label{eq:rbf_rk4_lte_expd}
        = {} & \left( 1 - b_1 - b_2 - b_3 - b_4 \right) f + \left( \A^4_1 f_t + \A^4_1 f_u f \right) h + \\
          {} & \Big[ \A^4_2 f_{tt} + 2 \A^4_2 f_{tu} f + \B^4_2 f_t f_u + \B^4_2 f^2_u f + \A^4_2 f_{uu} f^2 + \\
          {} & \left( b_2 c^2_2 \epsilon^2_{n2} + b_3 c^2_3 \epsilon^2_{n3} + b_4 c^2_4 \epsilon^2_{n4} \right) f_u u_n \Big] h^2+\\
          {} & \Big[ \A^4_3 f_{ttt} + 3 \A^4_3 f_{ttu} f + \B^4_3 f_t f_{tu} + 3 \A^4_3 f_{tuu} f^2 + \D^4_3 f_{tt} f_u + 
               \left( \B^4_3 + 2\D^4_3 \right) f_{tu} f_u f  + \\
          {} & \E^4_3 f_t f^2_u + \E^4_3 f^3_u f + \B^4_3 f_t f_{uu} f + \left( \B^4_3 + \D^4_3 \right) f_u f_{uu} f^2 + 
               \A^4_3 f_{uuu} f^3 + \\
          {} & \left( b_2 c^3_2 \epsilon^2_{n2} + b_3 c^3_3 \epsilon^2_{n3} + b_4 c^3_4 \epsilon^2_{n4} \right) \left( f_{tu} u_n + f_{uu} f u_n \right) + \\
          {} & \left( a_{32} b_3 c^2_2 \epsilon^2_{n2} + a_{42} b_4 c^2_2 \epsilon^2_{n2} + a_{43} b_4 c^2_3 \epsilon^2_{n3} \right) f^2_u u_n \Big] h^3 + \\
          {} & \bigg[ \A^4_4 f_{tttt} + 4 \A^4_4 f_{tttu} f + \B^4_4 f_t f_{ttu} + 6 \A^4_4 f_{ttuu} f^2 +  \D^4_4 f_{tt} f_{tu} +\\
          {} & 2 \D^4_4 f^2_{tu} f + 2 \B^4_4 f_t f_{tuu} f + 4 \A^4_4 f_{tuuu} f^3 + \E^4_4 f_{ttt} f_u + 
               \left( \B^4_4 + 3 \E^4_4 \right) f_{ttu} f_u f + \\  
          {} & \F^4_4 f_t f_{tu} f_u + \left( 2 \B^4_4 + 3 \E^4_4 \right) f_u f_{tuu} f^2 + \G^4_4 f_{tt} f^2_u + \\
          {} & \left( \F^4_4 + 2 \G^4_4 \right) f_{tu} f^2_u f + \HH^4_4 f_t^2 f_{uu} + \D^4_4 f_{tt} f_{uu} f + 3 \D^4_4 f_{tu} f_{uu} f^2 + \\
          {} & \left( \F^4_4 + 2 \HH^4_4 \right) f_t f_u f_{uu} f + \left( \F^4_4 + \G^4_4 + \HH^4_4 \right) f_u^2 f_{uu} f^2 + 
               \D^4_4 f^2_{uu} f^3 + \\
          {} & \B^4_4 f_t f_{uuu} f^2 + \left( \B^4_4 + \E^4_4 \right) f_u f_{uuu} f^3 + \A^4_4 f_{uuuu} f^4 + \\
          {} & \frac{1}{2} \left(b_2 c_2^4 \epsilon^2_{n2} + b_3 c_3^4 \epsilon^2_{n3} + b_4 c_4^4 \epsilon^2_{n4} \right)
               \left( f_{ttu} u_n + 2 f_{tuu} f u_n + f_{uuu} f^2 u_n \right) + \\
          {} & \left( a_{32} b_3 c_2^3 \epsilon^2_{n2} + a_{42} b_4 c_2^3 \epsilon^2_{n2} + a_{43} b_4 c_3^3 \epsilon^2_{n3} + a_{32} b_3 c_2^2 c_3 \epsilon^2_{n2} + a_{42} b_4 c_2^2 c_4 \epsilon^2_{n2}  + a_{43} b_4 c_3^2 c_4 \epsilon^2_{n3} \right) \left( f_{tu} f_u u_n + f_u f_{uu} f u_n \right) + \\
          {} & \left( a_{32} b_3 c_2 c_3^2 \epsilon^2_{n3} + a_{42} b_4 c_2 c_4^2 \epsilon^2_{n4} + a_{43} b_4 c_3 c_4^2 \epsilon^2_{n4} \right) \left( f_t f_{uu} u_n + f_u f_{uu} f u_n \right)  \\     
          {} & - \frac{1}{2} \left( b_2 c_2^4 \epsilon_{n2}^4 + b_3 c_3^4 \epsilon_{n3}^4 + b_4 c_4^4 \epsilon_{n4}^4 \right)  \left( f_u u_n + f_{uu} u_n^2 \right) + \\
          {} & \epsilon^2_{n2} a_{32} a_{43} b_4 c_2^2 f_u^3 u_n + \frac{1}{120} f^3_u u''_n \bigg] h^4 + O(h^5),   
 \end{split}
\end{align}
with $\A^4_1, \, \A^4_2, \, \B^4_2, \, \A^4_3, \, \B^4_3, \, \D^4_3, \, \E^4_3, \, \A^4_4, \, \B^4_4, \, \D^4_4, \, \E^4_4, \, \F^4_4, \, \G^4_4, \, \HH^4_4$ given in Appendix.

If we combine the conditions \eqref{eq:rk_consistent} and \eqref{eq:rk4_cond} with the extra restrictions in $h^2$ and $h^3$ terms
\begin{align}
\begin{split} \label{eq:rbf_rk4_cond}
 b_2 c^2_2 \epsilon^2_{n2} + b_3 c^2_3 \epsilon^2_{n3} + b_4 c^2_4 \epsilon^2_{n4} &= 0, \\
 b_2 c^3_2 \epsilon^2_{n2} + b_3 c^3_3 \epsilon^2_{n3} + b_4 c^3_4 \epsilon^2_{n4} &= 0, \\
 a_{32} b_3 c^2_2 \epsilon^2_{n2} + a_{42} b_4 c^2_2 \epsilon^2_{n2} + a_{43} b_4 c^2_3 \epsilon^2_{n3} &= 0,
\end{split}
\end{align}
then $\tau_n = O(h^4)$ and the RBF Runge-Kutta method \eqref{eq:rbf_rk4} remains fourth-order accurate.
The annihilation of the $h^4$ term is necessary to attain order $5$, and here we provide two possible solutions:
\begin{enumerate}[label=\Roman*.]
\item $a_{32} b_3 c_2 c_3^2 + a_{42} b_4 c_2 c_4^2 + a_{43} b_4 c_3 c_4^2 = \frac{1}{10}, \, a_{32} b_3 c_2^2 c_3 + a_{42} b_4 c_2^2 c_4 + a_{43} b_4 c_3^2 c_4 = \frac{1}{15}$ \\
      The coefficients of $f_t f_{ttu}, \, f_{tt} f_{tu}, \, f_{tu}^2 f, \, f_t f_{tuu} f, \, f_{tt} f_u^2, \, f_{tt} f_{uu} f, \, f_{tu} f_{uu} f^2, \, f_{uu}^2 f^3$ and $f_t f_{uuu} f^2$ in $h^4$ are zero.
      Then $a_{21} = c_2 = \frac{2}{5}, \, a_{31} = -\frac{3}{20}, \, a_{32} = \frac{3}{4}, \, c_3 = \frac{3}{5}, \, a_{41} = \frac{19}{44}, \, a_{42} = -\frac{15}{44}, \, a_{43} = \frac{10}{11}, \, c_4 = 1, \, b_1 = \frac{11}{72}, \, b_2 = \frac{25}{72}, \, b_3 = \frac{25}{72}$ and $b_4 = \frac{11}{72}$.
      By \eqref{eq:rbf_rk4_cond}, we have
$$
   \epsilon^2_{n3} = - \frac{2}{3} \epsilon^2_{n2}, \quad \epsilon^2_{n4} = \frac{2}{11} \epsilon^2_{n2}.
$$       
      Equating the rest of terms in $h^4$ to zero gives 
\begin{equation} \label{eq:qd_eqn_I}
 \alpha \epsilon^4_{n2} + \beta \epsilon^2_{n2} + \gamma = 0,
\end{equation} 
where 
\begin{align*}
 \alpha = {} & 672 \left( f_u + f_{uu} u_n \right) u_n, \\
 \beta  = {} & - \left( 132 f_{ttu} + 264 f_{tuu} f - 924 f_{tu} f_u - 540 f_t f_{uu} - 1464 f_u f_{uu} f + 132 f_{uuu} f^2 + 660 f_u^3 \right) u_n, \\
 \gamma = {} & 11 f_{tttt} + 44 f_{tttu} f + 66 f_{ttuu} f^2 + 44 f_{tuuu} f^3 - 44 f_{ttt} f_u - 132 f_{ttu} f_u f + \\
     {} & 330 f_t f_{tu} f_u - 132 f_u f_{tuu} f^2 + 330 f_{tu} f_u^2 f + 135 f_t^2 f_{uu} + 600 f_t f_u f_{uu} f + \\
     {} & 465 f_u^2 f_{uu} f^2 - 44 f_u f_{uuu} f^3 + 11 f_{uuuu} f^4 - 330 f_u^3 u''_n.
\end{align*}
\item $a_{32} b_3 c_2 c_3^2 + a_{42} b_4 c_2 c_4^2 + a_{43} b_4 c_3 c_4^2 = \frac{1}{10}, \, a_{32}^2 c_2^2 b_3 + \left( a_{42} c_2 + a_{43} c_3 \right)^2 b_4 = \frac{1}{20}$ \\
      The $f_t f_{ttu}, \, f_t f_{tuu} f, \, f_t^2 f_{uu}$ and $f_t f_{uuu} f^2$ terms in $h^4$ vanish.
      Then $a_{21} = c_2 = \frac{1}{4}, \, a_{31} = -\frac{6}{25}, \, a_{32} = \frac{21}{25}, \, c_3 = \frac{3}{5}, \, a_{41} = \frac{6}{5}, \, a_{42} = -\frac{57}{35}, \, a_{43} = \frac{10}{7}, \, c_4 = 1, \: b_1 = \frac{1}{9}, \, b_2 = \frac{16}{63}, \, b_3 = \frac{125}{252}$ and $b_4 = \frac{5}{36}$.
      Solving \eqref{eq:rbf_rk4_cond}, we find
$$
   \epsilon^2_{n3} = - \frac{1}{6} \epsilon^2_{n2}, \quad \epsilon^2_{n4} = \frac{1}{10} \epsilon^2_{n2}.
$$
	  If we set the leading error term in $O(h^4)$ to zero, then  
\begin{equation} \label{eq:qd_eqn_II}
 \alpha \epsilon^4_{n2} + \beta \epsilon^2_{n2} + \gamma =0, 
\end{equation}
where 
\begin{align*}
 \alpha = {} & 12 \left( f_u + f_{uu} u_n \right) u_n, \\
 \beta  = {} & - \left( 12 f_{ttu} + 24 f_{tuu} f - 84 f_{tu} f_u - 84 f_u f_{uu} f + 12 f_{uuu} f^2 + 60 f_u^3 \right) u_n, \\
 \gamma = {} & f_{tttt} + 4 f_{tttu} f + 6 f_{ttuu} f^2 + 18 f_{tt} f_{tu} + 36 f_{tu}^2 f + 4 f_{tuuu} f^3 - 4 f_{ttt} f_u - \\
          {} & 12 f_{ttu} f_u f + 48 f_t f_{tu} f_u - 12 f_u f_{tuu} f^2 - 18 f_{tt} f_u^2 + 12 f_{tu} f_u^2 f + 18 f_{tt}f_{uu}f+\\
          {} & 54 f_{tu} f_{uu}f^2 + 48 f_t f_u f_{uu} f + 30 f_u^2 f_{uu} f^2 + 18 f_{uu}^2 f^3 - 4 f_u f_{uuu} f^3 -48f_u^3 u''_n.
\end{align*}
\end{enumerate}
\begin{note}
 It is not always viable to vanish the $h^4$ term as the quadratic equation above cannot guarantee a real root, in which case the four-stage RBF Runge-Kutta method \eqref{eq:rbf_rk4} is still of order 4.
\end{note}

\section{Convergence of RBF Runge-Kutta methods} \label{sec:convg}
In this section, we prove the convergence of the proposed RBF Runge-Kutta methods in Sect. \ref{sec:rbf_rk}.

In order to guarantee that there is a unique solution to the initial value problem \eqref{eq:ivp}, some form of continuity is required in the function $f(t,u)$.
We say that the function $f(t,u)$ is Lipschitz continuous in $u$ if there exists some constant $L$ such that
$$
   \left| f(t, u_1) - f(t, u_2) \right| \leqslant L |u_1 - u_2|.
$$
Now we want to solve the initial value problem \eqref{eq:ivp} numerically on the interval $[a,b]$. 
This interval $[a,b]$ can divided by a uniform grid of $N+1$ points
$$
   a = t_0 < t_1 < \cdots < x_N = b,
$$
where the grid points are defined as $t_n = a + n h, \, n = 0, 1, \cdots, N$ and $h = (b-a)/N$ is the step size. 
Let $u_n = u(t_n)$ be the value of the exact solution and $v_n$ the numerical approximation by applying the RBF Runge-Kutta methods in succession.

\subsection{Two-stage methods}
\begin{theorem} \label{thm:rk2}
Suppose that $\epsilon^2_n$ is bounded for all $n = 0, 1, \cdots, N-1$.
Then the two-stage RBF Runge-Kutta method \eqref{eq:rbf_rk2}, satisfying the conditions \eqref{eq:rk_consistent}, \eqref{eq:rk_cond} and \eqref{eq:rk2_cond}, is convergent.
\end{theorem}

\begin{proof}
We fix $t = t_n$ and consider the error $E_{n+1}$ between $u_{n+1}$ and $v_{n+1}$.
Let 
\begin{align*}
 k_1(w) &= f(t_n, w), \\
 k_2(w) &= f(t_n + c_2 h, w e^{- \epsilon^2_n (c_2 h)^2} + h c_2 k_1).
\end{align*}
From the local truncation error \eqref{eq:rbf_rk2_lte}, the exact solution satisfies
\begin{equation} \label{eq:rbf_rk2_lte_exact}
 u_{n+1} = u_n + h \left( b_1 k_1(u_n) + b_2 k_2(u_n) \right) + h \tau_n. 
\end{equation}
By the Lipschitz continuity of $f(t,u)$ in $u$, we have
$$
   |k_1(u_n) - k_1(v_n)| = |f(t_n, u_n) - f(t_n, v_n)| \leqslant L |u_n - v_n|,
$$
and 
\begin{align*}
    |k_2(u_n) - k_2(v_n)| 
 &= \left| f \left( t_n + c_2 h, u_n e^{- \epsilon^2_n (c_2 h)^2} + h c_2 k_1(u_n) \right) - f \left( t_n + c_2 h, v_n e^{- \epsilon^2_n (c_2 h)^2} + h c_2 k_1(v_n) \right) \right| \\
 &\leqslant L \left| u_n e^{- \epsilon^2_n c_2^2 h^2} + c_2 h k_1(u_n) - v_n e^{- \epsilon^2_n c_2^2 h^2} - c_2 h k_1(v_n) \right|\\
 &\leqslant L \left( |u_n - v_n| e^{- \epsilon^2_n c_2^2 h^2} + c_2 h |k_1(u_n) - k_1(v_n)| \right) \\
 &\leqslant L \left( |u_n - v_n| e^{- \epsilon^2_n c_2^2 h^2} + c_2 h L |u_n - v_n| \right) \\
 &= L \left( e^{- \epsilon^2_n c_2^2 h^2} + c_2 h L \right) |u_n - v_n|.
\end{align*}
We subtract \eqref{eq:rbf_rk2_lte_exact} from \eqref{eq:rbf_rk2}, which gives
\begin{align}
 E_{n+1} &= |u_{n+1} - v_{n+1}| \nonumber \\
         &= \left| (u_n-v_n) + h \left[ b_1 (k_1(u_n)-k_1(v_n)) + b_2 (k_2(u_n)-k_2(v_n)) \right] + h \tau_n \right| \nonumber \\
         &\leqslant |u_n - v_n| + h \left[ b_1 |k_1(u_n)-k_1(v_n)| + b_2 |k_2(u_n)-k_2(v_n)| \right] + h |\tau_n| \nonumber \\
         &\leqslant |u_n - v_n| + h \left[ b_1 L |u_n - v_n| + b_2 L \left( e^{- \epsilon^2_n c_2^2 h^2} + c_2 h L \right) |u_n - v_n| \right] + h |\tau_n| \nonumber \\
         &= \varphi_n E_n + h |\tau_n|, \label{eq:rbf_rk2_ineq}
\end{align}
with $\varphi_n = 1 + b_1 h L + b_2 h L e^{- \epsilon^2_n c_2^2 h^2} + \frac{1}{2} h^2 L^2$ and $E_n = |u_n - v_n|$.
Let
$$
   \C_2 = \sup_{\substack{0<h\leqslant b-a\\n=0,\cdots,N-1}} \frac{b_2 L \left( e^{- \epsilon^2_n c_2^2 h^2} - 1 \right)}{e^{hL}}.
$$
Since $\epsilon^2_n$ is bounded for all $n = 0, 1, \cdots, N-1$, then $\C_2 < \infty$, and thus
\begin{align*}
 \varphi_n &= 1 + h L + \frac{1}{2} h^2 L^2 + b_2 h L \left( e^{- \epsilon^2_n c_2^2 h^2} - 1 \right) \\
           &\leqslant (1+\C_2 h) e^{hL} \\
           &\leqslant e^{h(L+\C_2)}.
\end{align*}
From the inequality \eqref{eq:rbf_rk2_ineq}, we obtain by induction
\begin{equation} \label{eq:rbf_rk2_errn}
 E_n \leqslant \prod_{j=0}^{n-1} \varphi_j E_0 + h \sum_{j=1}^{n-1} \left( \prod_{m=j}^{n-1} \varphi_m \right) |\tau_{j-1}| + h |\tau_{n-1}|.
\end{equation}
Note that
$$
   \prod_{m=j}^{n-1} \varphi_m \leqslant \prod_{m=0}^{n-1} \varphi_m \leqslant \left( e^{h(L+\C_2)} \right)^n \leqslant \left( e^{h(L+\C_2)} \right)^N = e^{(L+\C_2)(b-a)}.
$$
With the settings $E_0 = 0$ and 
$$ ||\tau||_{\infty} = \max_{n=0,\cdots,N-1} |\tau_n|, $$
it follows from \eqref{eq:rbf_rk2_errn} that for every $n = 0, 1, \cdots, N$,
$$
   E_n \leqslant h e^{(L+\C_2)T} \sum_{j=0}^{n-1} ||\tau||_{\infty} = n h e^{(L+\C_2)(b-a)} ||\tau||_{\infty} \leqslant (b-a) e^{(L+\C_2)(b-a)} ||\tau||_{\infty},
$$
and hence
$$
   \lim_{\substack{h \to 0\\Nh = (b-a)}} E_N = 0.
$$
Therefore, $v_N \to u_N = u(b)$ as $h \to 0$ and the two-stage RBF Runge-Kutta method converges. \qed
\end{proof}

\subsection{Three-stage methods}
\begin{theorem} \label{thm:rk3}
Suppose that $\epsilon^2_{n2}$ and $\epsilon^2_{n3}$ are bounded for all $n = 0, 1, \cdots, N-1$.
Then the three-stage RBF Runge-Kutta method \eqref{eq:rbf_rk3}, satisfying the conditions \eqref{eq:rk_consistent}, \eqref{eq:rk_cond} and \eqref{eq:rk3_cond}, converges.
\end{theorem}

\begin{proof}
Fixing $t = t_n$, we define
\begin{align*}
 k_1(w) &= f(t_n, w), \\
 k_2(w) &= f(t_n + c_2 h, w e^{- \epsilon^2_{n2} (c_2 h)^2} + h c_2 k_1), \\
 k_3(w) &= f(t_n + c_3 h, w e^{- \epsilon^2_{n3} (c_3 h)^2} + h (a_{31} k_1 + a_{32} k_2)).
\end{align*}
By the definition of local truncation error \eqref{eq:rbf_rk3_lte}, 
\begin{equation} \label{eq:rbf_rk3_lte_exact}
 u_{n+1} = u_n + h \left( b_1 k_1(u_n) + b_2 k_2(u_n) + b_3 k_3(u_n) \right) + h \tau_n .
\end{equation}
Since $f(t,u)$ is Lipschitz continuous in $u$, we can estimate
\begin{align*}
 |k_1(u_n) - k_1(v_n)| &\leqslant L |u_n - v_n|, \\
 |k_2(u_n) - k_2(v_n)| &\leqslant L \left( e^{- \epsilon^2_{n2} c^2_2 h^2} + c_2 h L \right) |u_n - v_n|,
\end{align*}
and 
\begin{align*}
 |k_3(u_n) - k_3(v_n)| &\leqslant L \left| u_n e^{- \epsilon^2_{n3} c^2_3 h^2} + a_{31} h k_1(u_n) + a_{32} h k_2(u_n) - v_n e^{- \epsilon^2_{n3} c^2_3 h^2} - a_{31} h k_1(v_n) - a_{32} h k_2(v_n) \right| \\
                       &\leqslant L \left[ |u_n - v_n| e^{- \epsilon^2_{n3} c^2_3 h^2} + a_{31} h |k_1(u_n) - k_1(v_n)| + a_{32} h |k_2(u_n) - k_2(v_n)| \right] \\
                       &\leqslant L \left[ e^{- \epsilon^2_{n3} c^2_3 h^2} + \left( a_{31} + a_{32} e^{- \epsilon^2_{n2} c^2_2 h^2}  \right) h L + a_{32} c_2 h^2 L^2 \right] |u_n - v_n|.
\end{align*}
Subtracting \eqref{eq:rbf_rk3_lte_exact} from \eqref{eq:rbf_rk3} and using the inequalities above yield
\begin{equation} \label{eq:rbf_rk3_ineq}
 E_{n+1} = |u_{n+1} - v_{n+1}| \leqslant \varphi_n E_n + h |\tau_n|,
\end{equation}
where $E_n = |u_n - v_n|$ and
$$
 \varphi_n = 1 + \left( b_1 + b_2 e^{- \epsilon^2_{n2} c^2_2 h^2} + b_3 e^{- \epsilon^2_{n3} c^2_3 h^2} \right) h L 
               + \left( b_2 c_2 + a_{31} b_3 + a_{32} b_3 e^{- \epsilon^2_{n2} c^2_2 h^2} \right) h^2 L^2 + \frac{1}{6} h^3 L^3.
$$
Let
$$
   \C_3 = \sup_{\substack{0<h\leqslant b-a\\n=0,\cdots,N-1}} \frac{(b_2 L + a_{32} b_3 h L^2) \left( e^{-\epsilon^2_{n2} c^2_2 h^2} - 1 \right) + b_3 L \left( e^{- \epsilon^2_{n3} c^2_3 h^2} - 1 \right)}{e^{hL}}.
$$
Then $\C_3$ is bounded since $\epsilon^2_{n2}$ and $\epsilon^2_{n3}$ are bounded for all $n = 0, 1, \cdots, N-1$.
By the conditions \eqref{eq:rk_consistent}, \eqref{eq:rk_cond} and \eqref{eq:rk3_cond}, 
\begin{align*}
 \varphi_n &= 1 + h L + \frac{1}{2} h^2 L^2 + \frac{1}{6} h^3 L^3 + (b_2 h L + a_{32} b_3 h^2 L^2) \left( e^{-\epsilon^2_{n2} c^2_2 h^2} - 1 \right) + b_3 h L \left( e^{- \epsilon^2_{n3} c^2_3 h^2} - 1 \right) \\
           &\leqslant (1+\C_3 h) e^{hL} \\
           &\leqslant e^{h(L+\C_3)}.
\end{align*}
The rest of the argument is similar to the proof of Theorem \ref{thm:rk2} by replacing $\C_2$ with $\C_3$.
Therefore, we conclude $v_N \to u_N = u(b)$ as $h \to 0$ and the three-stage RBF Runge-Kutta method converges. \qed
\end{proof}

\subsection{Four-stage RBF Runge-Kutta methods}
\begin{theorem} \label{thm:rk4}
Suppose that $\epsilon^2_{n2}, \, \epsilon^2_{n3}$ and $\epsilon^2_{n4}$ are bounded for all $n = 0, 1, \cdots, N-1$.
Then the four-stage RBF Runge-Kutta method \eqref{eq:rbf_rk4}, satisfying the conditions \eqref{eq:rk_consistent}, \eqref{eq:rk_cond} and \eqref{eq:rk4_cond}, converges to the exact solution.
\end{theorem}

\begin{proof}
For $t = t_n$, define
\begin{align*}
 k_1(w) &= f(t_n, w), \\
 k_2(w) &= f(t_n + c_2 h, w e^{- \epsilon^2_{n2} (c_2 h)^2} + h c_2 k_1), \\
 k_3(w) &= f(t_n + c_3 h, w e^{- \epsilon^2_{n3} (c_3 h)^2} + h (a_{31} k_1 + a_{32} k_2)), \\
 k_4(w) &= f(t_n + c_4 h, w e^{- \epsilon^2_{n4} (c_4 h)^2} + h (a_{41} k_1 + a_{42} k_2 + a_{43} k_3)). 
\end{align*}
We can write the exact solution as
\begin{equation} \label{eq:rbf_rk4_lte_exact}
 u_{n+1} = u_n + h \left( b_1 k_1(u_n) + b_2 k_2(u_n) + b_3 k_3(u_n) + b_4 k_4(u_n) \right) + h \tau_n,
\end{equation}
from the local truncation error \eqref{eq:rbf_rk4_lte}. 
The Lipschitz continuity of $f(t,u)$ in $u$ gives
\begin{align*}
 |k_1(u_n) - k_1(v_n)| &\leqslant L |u_n - v_n|, \\
 |k_2(u_n) - k_2(v_n)| &\leqslant L \left( e^{- \epsilon^2_{n2} c^2_2 h^2} + c_2 h L \right) |u_n - v_n|, \\
 |k_3(u_n) - k_3(v_n)| &\leqslant L \left[ e^{- \epsilon^2_{n3} c^2_3 h^2} + \left( a_{31} + a_{32} e^{- \epsilon^2_{n2} c^2_2 h^2}  \right) h L + a_{32} c_2 h^2 L^2 \right] |u_n - v_n|,
\end{align*}
and 
\begin{align*}
 |k_4(u_n) - k_4(v_n)| &\leqslant L \left[ |u_n - v_n| e^{- \epsilon^2_{n4} c^2_4 h^2} + a_{41} h |k_1(u_n) - k_1(v_n)| + a_{42} h |k_2(u_n) - k_2(v_n)| + a_{43} h |k_3(u_n) - k_3(v_n)| \right] \\
                       &\leqslant L \left[ e^{- \epsilon^2_{n4} c^2_4 h^2} + \left( a_{41} + a_{42} e^{- \epsilon^2_{n2} c^2_2 h^2} + a_{43} e^{- \epsilon^2_{n3} c^2_3 h^2} \right) h L \right.\\
                       & ~~~~~~~\left. + \left( a_{42} c_2 + a_{31} a_{43} + a_{32} a_{43} e^{- \epsilon^2_{n2} c^2_2 h^2} \right) h^2 L^2 + a_{32} a_{43} c_2 h^3 L^3 \right] |u_n - v_n|.
\end{align*}
By subtracting the two expressions \eqref{eq:rbf_rk4_lte_exact} and \eqref{eq:rbf_rk4} and using the inequalities above, we obtain
\begin{equation} \label{eq:rbf_rk4_ineq}
 E_{n+1} = |u_{n+1} - v_{n+1}| \leqslant \varphi_n E_n + h |\tau_n|,
\end{equation}
where $E_n = |u_n - v_n|$ and
\begin{align*}
 \varphi_n = 1 &+ \left( b_1 + b_2 e^{- \epsilon^2_{n2} c^2_2 h^2} + b_3 e^{- \epsilon^2_{n3} c^2_3 h^2} + b_4 e^{- \epsilon^2_{n4}c^2_4 h^2} \right) h L \\
               &+ \left( b_2 c_2 + a_{31} b_3 + a_{32} b_3 e^{- \epsilon^2_{n2} c^2_2 h^2}) + a_{41} b_4 + a_{42} b_4 e^{-\epsilon^2_{n2} c^2_2 h^2} + a_{43} b_4 e^{- \epsilon^2_{n3} c^2_3 h^2} \right) h^2 L^2 \\
               &+ \left( a_{32} b_3 c_2 + a_{42} b_4 c_2 + a_{31} a_{43} b_4 + a_{32} a_{43} b_4 e^{- \epsilon^2_{n2} c^2_2 h^2} \right) h^3 L^3 + \frac{1}{24} h^4 L^4.
\end{align*} 
Let
\begin{align*}
 \C_4 = \sup_{\substack{0<h\leqslant b-a\\n=0,\cdots,N-1}} \frac{1}{e^{hL}} & \left\{ [b_2 L + (a_{32} b_3 + a_{42} b_4) h L^2 + a_{32} a_{43} b_4 h^2 L^3] \left( e^{-\epsilon^2_{n2} c^2_2 h^2} - 1 \right) \right. \\
 & \; \left. + (b_3 L + a_{43} b_4 h L^2) \left( e^{- \epsilon^2_{n3} c^2_3 h^2} - 1 \right) + b_4 L \left( e^{- \epsilon^2_{n4} c^2_4 h^2} - 1 \right) \right\}.
\end{align*}   
Then $\C_4$ is bounded assuming that $\epsilon^2_{n2}, \, \epsilon^2_{n3}$ and $\epsilon^2_{n4}$ are bounded for all $n$.
According to the conditions \eqref{eq:rk_consistent}, \eqref{eq:rk_cond} and \eqref{eq:rk4_cond}, we have
$$
   \varphi_n \leqslant (1+\C_4 h) e^{hL} \leqslant e^{h(L+\C_4)}.
$$
The rest of the proof then follows as in the proof of Theorem \ref{thm:rk2} with $\C_4$ in place of $\C_2$.
Therefore, $v_N \to u_N = u(b)$ as $h \to 0$, which implies that the four-stage RBF Runge-Kutta method converges. \qed
\end{proof}

\section{Stability regions} \label{sec:sregion}
We have constructed the RBF Runge-Kutta methods up to four stages and proved the convergence of those methods over some bounded interval as the step size $h$ tends to zero. 
However, the convergence of a numerical method cannot guarantee the reasonable results with a given value $h>0$ in practice, which is related to the stability region. 
In the present section, we plot the stability region of each RBF Runge-Kutta method and compare with the one of its correspondent Runge-Kutta method.

All Runge-Kutta methods presented in this paper are one-step methods.
Applying a one-step method to the test problem $u' = \lambda u$ yields an expression of form $v_{n+1} = R(z) v_n$, where $R(z)$ is the stability function of $z = \lambda h$.
If the one-step method is consistent, then $R(z)$ approximates $e^z$ around $z=0$.
If the one-step method has order $p$, then
$$
   R(z) - e^z = O(z^{p+1}) \text{ as } z \to 0.
$$
The stability functions of the two-, three- and four-stage Runge-Kutta methods are 
\begin{align*}
 R(z) &= 1 + z + \frac{1}{2} z^2, \\
 R(z) &= 1 + z + \frac{1}{2} z^2 + \frac{1}{6} z^3, \\
 R(z) &= 1 + z + \frac{1}{2} z^2 + \frac{1}{6} z^3 + \frac{1}{24} z^4,
\end{align*}
respectively.
The stability function of the two-stage RBF Runge-Kutta method in Subsect. \ref{sec:rbf_rk2} is
$$ R(z) = 1 + \left( \frac{1}{4} + \frac{3}{4} e^{\frac{2}{9} z^2} \right) z + \frac{1}{2} z^2. $$
The stability functions of the three-stage RBF Runge-Kutta methods in Subsect. \ref{sec:rbf_rk3} are
\begin{enumerate}[label=\Roman*.]
\item $$ R(z) = 1 + \left( \frac{1}{6} + \frac{2}{3} e^{\frac{1}{8} z^2} + \frac{1}{6} e^{-\frac{1}{2} z^2} \right) z + \left( \frac{1}{6} + \frac{1}{3} e^{\frac{1}{8} z^2} \right) z^2 + \frac{1}{6} z^3. $$
\item 
 \begin{enumerate}
 \item \begin{align*} 
        R(z) = 1 + {} & \left( \frac{1}{8} + \frac{77+3 \sqrt{33}}{176} e^{\frac{15-\sqrt{33}}{96} z^2} + \frac{77-3 \sqrt{33}}{176} e^{-\frac{111+\sqrt{33}}{768} z^2} \right) z + \\
                   {} & \left( \frac{9-\sqrt{33}}{48} + \frac{15+\sqrt{33}}{48} e^{\frac{15-\sqrt{33}}{96} z^2} \right) z^2 + \frac{1}{6} z^3.
       \end{align*} 
 \item \begin{align*} 
        R(z) = 1 + {} & \left( \frac{1}{8} + \frac{77-3 \sqrt{33}}{176} e^{\frac{15+\sqrt{33}}{96} z^2} + \frac{77+3 \sqrt{33}}{176} e^{\frac{-111+\sqrt{33}}{768} z^2} \right) z + \\
                   {} & \left( \frac{9+\sqrt{33}}{48} + \frac{15-\sqrt{33}}{48} e^{\frac{15+\sqrt{33}}{96} z^2} \right) z^2 + \frac{1}{6} z^3. 
       \end{align*}
 \end{enumerate}
\item
 \begin{enumerate}
 \item $$ R(z) = 1 + \left( \frac{1}{10} + \frac{1}{2} e^{{1}{12} z^2} + \frac{2}{5} e^{-\frac{5}{48} z^2} \right) z + \frac{1}{2} e^{\frac{1}{12} z^2} z^2 + \frac{1}{6} z^3. $$  
 \item $$ R(z) = 1 + \left( \frac{1}{6} + \frac{1}{6} e^{\frac{1}{4} z^2} + \frac{2}{3} e^{-\frac{1}{16} z^2} \right) z + \left( \frac{1}{3} + \frac{1}{6} e^{\frac{1}{4} z^2} \right) z^2 + \frac{1}{6} z^3. $$ 
 \end{enumerate}
\item $$ R(z) = 1 + \left( \frac{2}{9} + \frac{1}{3} e^{\frac{1}{8} z^2} + \frac{4}{9} e^{-\frac{3}{32} z^2} \right) z + \left( \frac{1}{6} + \frac{1}{3} e^{\frac{1}{8} z^2} \right) z^2 + \frac{1}{6} z^3. $$
\end{enumerate}
The stability functions of the four-stage RBF Runge-Kutta methods in Subsect. \ref{sec:rbf_rk4} are
\begin{enumerate}[label=\Roman*.]
\item ($+$)
  \begin{align*}
   R(z) = 1 + {} & \left( \frac{11}{72} + \frac{25}{72} e^{-\frac{55+\sqrt{9185}}{700} z^2} + \frac{25}{72} e^{\frac{165+3\sqrt{9185}}{1400} z^2} + \frac{11}{72} e^{-\frac{55+\sqrt{9185}}{616} z^2} \right) z + \\
              {} & \left( \frac{11}{72} + \frac{25}{96} e^{-\frac{55+\sqrt{9185}}{700} z^2} - \frac{5}{96} e^{-\frac{55+\sqrt{9185}}{700} z^2} + \frac{5}{36} e^{\frac{165+3\sqrt{9185}}{1400} z^2} \right) z^2 + \\
              {} & \left( \frac{1}{16} + \frac{5}{48} e^{-\frac{55+\sqrt{9185}}{700} z^2} \right) z^3 + \frac{1}{24} z^4. 
  \end{align*}
      ($-$)
  \begin{align*}
   R(z) = 1 + {} & \left( \frac{11}{72} + \frac{25}{72} e^{\frac{-55+\sqrt{9185}}{700} z^2} + \frac{25}{72} e^{\frac{165-3\sqrt{9185}}{1400} z^2} + \frac{11}{72} e^{\frac{-55+\sqrt{9185}}{616} z^2} \right) z + \\
              {} & \left( \frac{11}{72} + \frac{25}{96} e^{\frac{-55+\sqrt{9185}}{700} z^2} - \frac{5}{96} e^{\frac{-55+\sqrt{9185}}{700} z^2} + \frac{5}{36} e^{\frac{165-3\sqrt{9185}}{1400} z^2} \right) z^2 + \\
              {} & \left( \frac{1}{16} + \frac{5}{48} e^{\frac{-55+\sqrt{9185}}{700} z^2} \right) z^3 + \frac{1}{24} z^4. 
  \end{align*}
\item ($+$)
	  \begin{align*}
       R(z) = 1 + {} & \left( \frac{1}{9} + \frac{16}{63} e^{-\frac{5+\sqrt{41}}{32} z^2} + \frac{125}{252} e^{\frac{15+3\sqrt{41}}{100} z^2} + \frac{5}{36} e^{-\frac{5+\sqrt{41}}{20} z^2} \right) z + \\
                  {} & \left( \frac{1}{9} + \frac{5}{12} e^{-\frac{5+\sqrt{41}}{32} z^2} - \frac{19}{84} e^{-\frac{5+\sqrt{41}}{32} z^2} + \frac{25}{126} e^{\frac{15+3\sqrt{41}}{100} z^2} \right) z^2 + \\
                  {} & \frac{1}{6} e^{-\frac{5+\sqrt{41}}{32} z^2} z^3 + \frac{1}{24} z^4. 
      \end{align*}
      ($-$)
      \begin{align*}
       R(z) = 1 + {} & \left( \frac{1}{9} + \frac{16}{63} e^{\frac{-5+\sqrt{41}}{32} z^2} + \frac{125}{252} e^{\frac{15-3\sqrt{41}}{100} z^2} + \frac{5}{36} e^{\frac{-5+\sqrt{41}}{20} z^2} \right) z + \\
                  {} & \left( \frac{1}{9} + \frac{5}{12} e^{\frac{-5+\sqrt{41}}{32} z^2} - \frac{19}{84} e^{\frac{-5+\sqrt{41}}{32} z^2} + \frac{25}{126} e^{\frac{15-3\sqrt{41}}{100} z^2} \right) z^2 + \\
                  {} & \frac{1}{6} e^{\frac{-5+\sqrt{41}}{32} z^2} z^3 + \frac{1}{24} z^4. 
      \end{align*} 
\end{enumerate}
The sign $\pm$ represents the solution of the quadratic equations \eqref{eq:qd_eqn_I} and \eqref{eq:qd_eqn_II}.

\begin{figure}[htbp]
 \centering
 \includegraphics[width=0.75\textwidth]{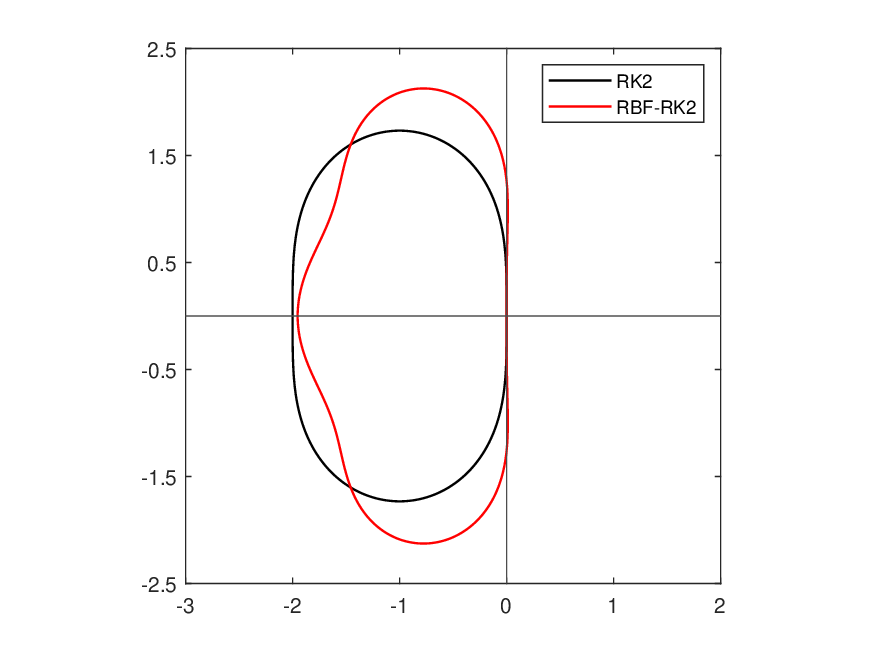}
 \caption{Stability regions of RK2 (black) and RBF-RK2 (red).}
 \label{fig:rk2}
\end{figure}

\begin{figure}[htbp]
 \centering
 \includegraphics[width=0.75\textwidth]{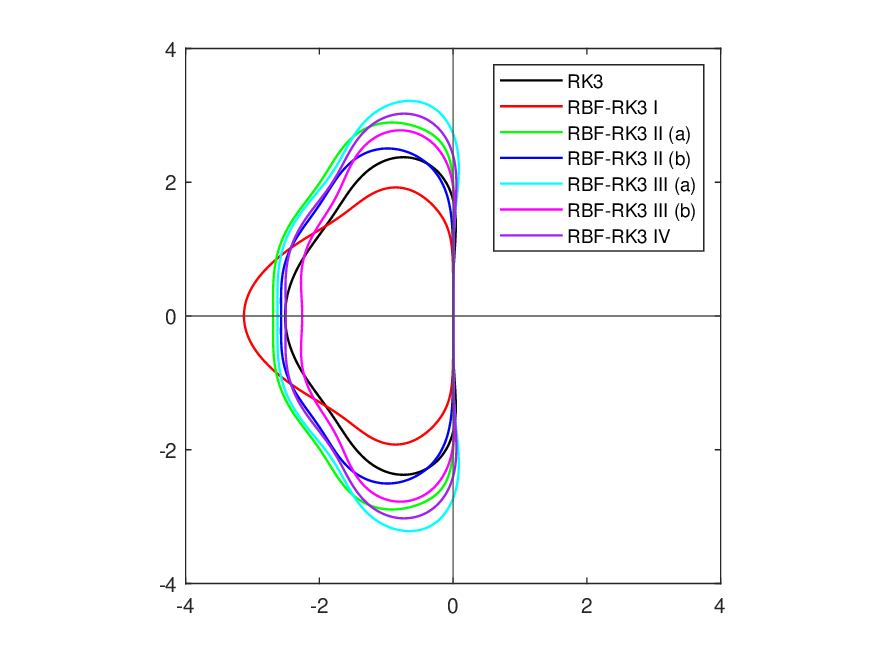}
 \caption{Stability regions of RK3 (black) and RBF-RK3 (red, green, blue, cyan, magenta and purple).}
 \label{fig:rk3}
\end{figure}

\begin{figure}[htbp]
 \centering
 \includegraphics[width=0.75\textwidth]{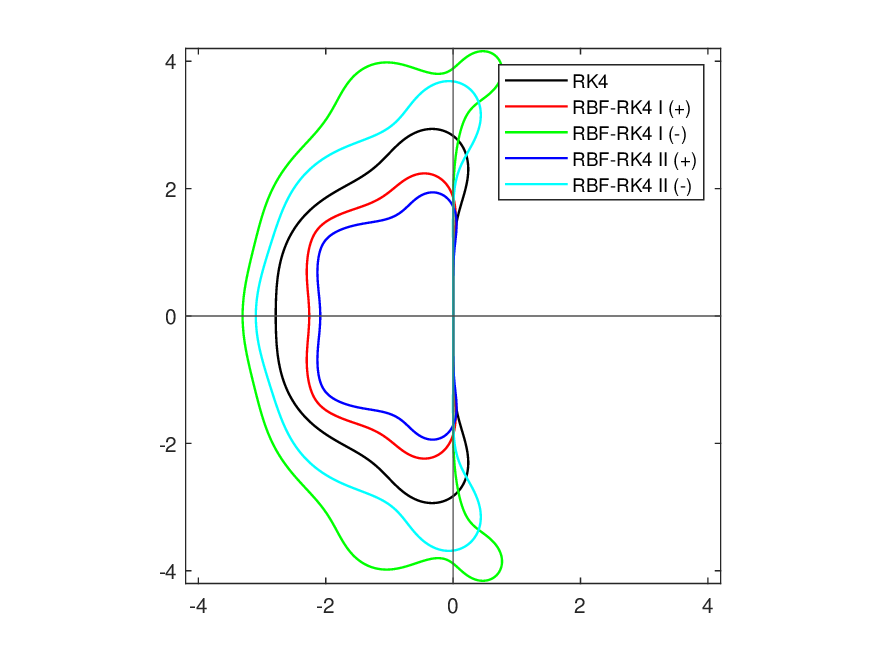}
 \caption{Stability regions of RK4 (black) and RBF-RK4 (red, green, blue and cyan).}
 \label{fig:rk4}
\end{figure}

Fig. \ref{fig:rk2} illustrates the stability regions for the two-stage Runge-Kutta method (RK2) and the two-stage RBF Runge-Kutta method (RBF-RK2).
The stability interval, which corresponds to the interval contained in both negative real axis and stability region, of RBF-RK2 is slightly smaller than RK2.
Fig. \ref{fig:rk3} plots the stability regions for the three-stage Runge-Kutta method (RK3) and the three-stage RBF Runge-Kutta methods (RBF-RK3).
From large to small interval of stability, the regions correspond to RBF-RK3 I, RBF-RK3 II (a), RBF-RK3 III (a), RBF-RK3 II (b), RBF-RK3 IV, RK3 and RBF-RK3 III (b).
Fig. \ref{fig:rk4} shows the stability regions for the four-stage Runge-Kutta method (RK4) and the four-stage RBF Runge-Kutta methods (RBF-RK4).
From large to small interval of stability, the regions correspond to RBF-RK4 I (-), RBF-RK4 II (-), RK4, RBF-RK4 I (+), and RBF-RK4 II (+).

\section{Numerical experiments} \label{sec:nr}
In this section we present some numerical experiments to show the performance and check the order of RBF Runge-Kutta methods introduced in Sect. \ref{sec:rbf_rk}, where we compare with the standard Runge-Kutta methods with all shape parameters equal to zero.

\begin{example} \label{ex:nusquare}
As the first example we consider
\begin{align*}
 \frac{du}{dt} &= -u^2, \quad t \geqslant 0, \\
          u(0) &= 1,
\end{align*}
where $f$ is only a function of $u$.
The exact solution is $u(t) = \frac{1}{t+1}$. 
We compute the numerical solution up to the final time $T=1$.

The global error at the final time versus $N$, as well as the order of accuracy, for the standard and RBF Runge-Kutta methods is displayed in Tables \ref{tab:nusquare_rk2}, \ref{tab:nusquare_rk3} and \ref{tab:nusquare_rk4}.
As can been seen, the $s$-stage Runge-Kutta method is $s$th-order accurate and the $s$-stage RBF Runge-Kutta method attains order $s+1$.
We also find that the global error of the $s$-stage RBF Runge-Kutta method is smaller than the one of the corresponding $s$-stage Runge-Kutta method at the same number of steps, except RBF-RK4 II (+) in Table \ref{tab:nusquare_rk4}.
However, the global error of the $s$-stage RBF Runge-Kutta method is typically larger than the one of the $(s\! +\! 1)$-stage Runge-Kutta method for the same $N$, though they achieve the same order of accuracy.
\end{example}

\begin{table}[h!]
\renewcommand{\arraystretch}{1.1}
\scriptsize
\centering
\caption{Global error and order of accuracy for Example \ref{ex:nusquare} by RK2 and RBF-RK2.}      
\begin{tabular}{clcrlc} 
\hline\noalign{\smallskip}
N & \multicolumn{2}{l}{RK2} & & \multicolumn{2}{l}{RBF-RK2} \\ 
    \cline{2-3}                 \cline{5-6}
  & Error & Order           & & Error & Order \\
\noalign{\smallskip}\hline\noalign{\smallskip}
10  & 9.34e-2 & --     & & 6.20e-5 & --     \\  
20  & 2.20e-4 & 2.0828 & & 7.10e-6 & 3.1257 \\  
40  & 5.36e-5 & 2.0410 & & 8.50e-7 & 3.0628 \\
80  & 1.32e-5 & 2.0204 & & 1.04e-7 & 3.0314 \\ 
160 & 3.28e-6 & 2.0102 & & 1.29e-8 & 3.0157 \\ 
320 & 8.17e-7 & 2.0051 & & 1.60e-9 & 3.0078 \\  
\noalign{\smallskip}\hline
\end{tabular}
\label{tab:nusquare_rk2}
\end{table}

\begin{table}[h!]
\renewcommand{\arraystretch}{1.1}
\scriptsize
\centering
\caption{Global error and order of accuracy for Example \ref{ex:nusquare} by RK3 and RBF-RK3.}      
\begin{tabular}{clcrlcrlcrlc} 
\hline\noalign{\smallskip}
N & \multicolumn{2}{l}{RK3 I} & & \multicolumn{2}{l}{RBF-RK3 I} & & \multicolumn{2}{l}{RK3 II (a)} & & \multicolumn{2}{l}{RBF-RK3 II (a)} \\
    \cline{2-3}                   \cline{5-6}                       \cline{8-9}                        \cline{11-12}  
  & Error & Order             & & Error & Order                 & & Error & Order                  & & Error & Order \\
\noalign{\smallskip}\hline\noalign{\smallskip}
10  & 1.93e-5  & --     & & 8.75e-7  & --     & & 3.14e-5  & --     & & 1.02e-6  & --     \\  
20  & 2.16e-6  & 3.1605 & & 4.58e-8  & 4.2573 & & 3.68e-6  & 3.0905 & & 6.16e-8  & 4.0496 \\  
40  & 2.57e-7  & 3.0752 & & 2.61e-9  & 4.1330 & & 4.46e-7  & 3.0450 & & 3.77e-9  & 4.0287 \\
80  & 3.13e-8  & 3.0363 & & 1.56e-10 & 4.0677 & & 5.49e-8  & 3.0225 & & 2.33e-10 & 4.0153 \\ 
160 & 3.86e-9  & 3.0178 & & 9.49e-12 & 4.0341 & & 6.81e-9  & 3.0112 & & 1.45e-11 & 4.0079 \\ 
320 & 4.80e-10 & 3.0088 & & 5.86e-13 & 4.0177 & & 8.48e-10 & 3.0056 & & 9.04e-13 & 4.0049 \\
\noalign{\smallskip}\hline\noalign{\smallskip}
N & \multicolumn{2}{l}{RK3 II (b)} & & \multicolumn{2}{l}{RBF-RK3 II (b)} & & \multicolumn{2}{l}{RK3 III (a)} & & \multicolumn{2}{l}{RBF-RK3 III (a)} \\
    \cline{2-3}                        \cline{5-6}                            \cline{8-9}                         \cline{11-12}  
  & Error & Order                  & & Error & Order                      & & Error & Order                   & & Error & Order \\
\noalign{\smallskip}\hline\noalign{\smallskip}
10  & 4.97e-5 & --     & & 2.30e-6  & --     & & 3.54e-5  & --     & & 1.53e-6  & -- \\  
20  & 5.76e-6 & 3.1107 & & 1.32e-7  & 4.1226 & & 4.16e-6  & 3.0897 & & 9.00e-8  & 4.0876 \\  
40  & 6.93e-7 & 3.0558 & & 7.91e-9  & 4.0627 & & 5.04e-7  & 3.0450 & & 5.45e-9  & 4.0459 \\
80  & 8.49e-8 & 3.0280 & & 4.84e-10 & 4.0317 & & 6.20e-8  & 3.0226 & & 3.35e-10 & 4.0235 \\ 
160 & 1.05e-8 & 3.0140 & & 2.99e-11 & 4.0159 & & 7.69e-9  & 3.0113 & & 2.08e-11 & 4.0118 \\ 
320 & 1.31e-9 & 3.0070 & & 1.86e-12 & 4.0081 & & 9.57e-10 & 3.0056 & & 1.29e-12 & 4.0059 \\
\noalign{\smallskip}\hline\noalign{\smallskip}
N & \multicolumn{2}{l}{RK3 III (b)} & & \multicolumn{2}{l}{RBF-RK3 III (b)} & & \multicolumn{2}{l}{RK3 IV} & & \multicolumn{2}{l}{RBF-RK3 IV} \\
    \cline{2-3}                         \cline{5-6}                             \cline{8-9}                    \cline{11-12}  
  & Error & Order                   & & Error & Order                       & & Error & Order              & & Error & Order \\
\noalign{\smallskip}\hline\noalign{\smallskip}
10  & 3.50e-5  & --     & & 2.30e-6  & --     & & 3.54e-5  & --     & & 1.65e-6  & -- \\  
20  & 4.14e-6  & 3.0794 & & 1.32e-7  & 4.1211 & & 4.16e-6  & 3.0899 & & 9.62e-8  & 4.1006 \\  
40  & 5.03e-7  & 3.0410 & & 7.93e-9  & 4.0617 & & 5.04e-7  & 3.0453 & & 5.80e-9  & 4.0518 \\
80  & 6.19e-8  & 3.0208 & & 4.85e-10 & 4.0311 & & 6.20e-8  & 3.0227 & & 3.56e-10 & 4.0262 \\ 
160 & 7.69e-9  & 3.0105 & & 3.00e-11 & 4.0156 & & 7.69e-9  & 3.0114 & & 2.21e-11 & 4.0132 \\ 
320 & 9.57e-10 & 3.0052 & & 1.86e-12 & 4.0075 & & 9.57e-10 & 3.0057 & & 1.37e-12 & 4.0063 \\
\noalign{\smallskip}\hline
\end{tabular}
\label{tab:nusquare_rk3}
\end{table}

\begin{table}[h!]
\renewcommand{\arraystretch}{1.1}
\scriptsize
\centering
\caption{Global error and order of accuracy for Example \ref{ex:nusquare} by RK4 and RBF-RK4.}      
\begin{tabular}{clcrlcrlc} 
\hline\noalign{\smallskip}
N & \multicolumn{2}{l}{RK4 I} & & \multicolumn{2}{l}{RBF-RK4 I (+)} & & \multicolumn{2}{l}{RBF-RK4 I (-)} \\
    \cline{2-3}                   \cline{5-6}                           \cline{8-9}
  & Error & Order             & & Error & Order                     & & Error & Order \\
\noalign{\smallskip}\hline\noalign{\smallskip}
10  & 2.44e-7  & --     & & 2.37e-7  & --     & & 4.51e-8  & --     \\  
20  & 1.69e-8  & 3.8546 & & 6.39e-9  & 5.2097 & & 1.30e-9  & 5.1115 \\   
40  & 1.09e-9  & 3.9501 & & 1.86e-10 & 5.1048 & & 3.92e-11 & 5.0569 \\ 
80  & 6.93e-11 & 3.9796 & & 5.60e-12 & 5.0523 & & 1.20e-12 & 5.0286 \\  
160 & 4.36e-12 & 3.9909 & & 1.72e-13 & 5.0285 & & 3.72e-14 & 5.0096 \\  
320 & 2.73e-13 & 3.9951 & & 5.77e-15 & 4.8930 & & 1.11e-15 & 5.0682 \\   
\noalign{\smallskip}\hline\noalign{\smallskip}
N & \multicolumn{2}{l}{RK4 II} & & \multicolumn{2}{l}{RBF-RK4 II (+)} & & \multicolumn{2}{l}{RBF-RK4 II (-)} \\
    \cline{2-3}                    \cline{5-6}                            \cline{8-9}
  & Error & Order              & & Error & Order                      & & Error & Order \\
\noalign{\smallskip}\hline\noalign{\smallskip}
10  & 6.13e-7  & --     & & 8.20e-7  & --     & & 5.55e-8  & --     \\  
20  & 3.74e-8  & 4.0352 & & 2.08e-8  & 5.3013 & & 1.58e-9  & 5.1392 \\   
40  & 2.30e-9  & 4.0237 & & 5.86e-10 & 5.1505 & & 4.69e-11  & 5.0694 \\ 
80  & 1.42e-10 & 4.0132 & & 1.74e-11 & 5.0751 & & 1.43e-12  & 5.0347 \\  
160 & 8.85e-12 & 4.0069 & & 5.29e-13 & 5.0373 & & 4.42e-14  & 5.0160 \\  
320 & 5.52e-13 & 4.0038 & & 1.63e-14 & 5.0234 & & 1.44e-15  & 4.9380 \\   
\noalign{\smallskip}\hline
\end{tabular}
\label{tab:nusquare_rk4}
\end{table}

\begin{example} \label{ex:stiff}
Next we proceed to solve the stiff problem
\begin{align*}
 \frac{du}{dt} &= -4 t^3 u^2, \quad t \geqslant -10, \\
        u(-10) &= 1/10001,
\end{align*}
with the exact solution $u(t) = 1/(t^4+1)$. 
The numerical solution will be finally approximated at the time $T=0$.

Tables \ref{tab:stiff_rk2}, \ref{tab:stiff_rk3} and \ref{tab:stiff_rk4} show the global error at $T=0$ for various values of $N$ and the order of accuracy.
All $s$-stage Runge-Kutta methods achieve $s$th-order accuracy and most of the $s$-stage RBF Runge-Kutta methods acquire $(s\! +\! 1)$th-order accuracy. 
However, RBF-RK3 III (b) in Table \ref{tab:stiff_rk3} shows the order increase, while RBF-RK4 II (-) in Table \ref{tab:stiff_rk4} exhibits an order reduction in this problem.
Similar to Example \ref{ex:nusquare}, the $s$-stage RBF Runge-Kutta method performs better than the corresponding $s$-stage Runge-Kutta method in terms of accuracy.
Moreover, the $s$-stage RBF Runge-Kutta method yields better accuracy than most of the $(s\! +\! 1)$-stage Runge-Kutta methods, which differs from Example \ref{ex:nusquare}.
\end{example}

\begin{table}[h!]
\renewcommand{\arraystretch}{1.1}
\scriptsize
\centering
\caption{Global error and order of accuracy for Example \ref{ex:stiff} by RK2 and RBF-RK2.}      
\begin{tabular}{clcrlc} 
\hline\noalign{\smallskip}
N & \multicolumn{2}{l}{RK2} & & \multicolumn{2}{l}{RBF RK2} \\ 
    \cline{2-3}                 \cline{5-6}
  & Error & Order           & & Error & Order \\
\noalign{\smallskip}\hline\noalign{\smallskip}
200  & 7.51e-1 & --     & & 3.56e-2 & --     \\  
400  & 4.40e-1 & 0.7726 & & 4.77e-3 & 2.8981 \\  
800  & 1.66e-1 & 1.4016 & & 6.11e-4 & 2.9663 \\
1600 & 4.79e-2 & 1.7972 & & 7.71e-5 & 2.9854 \\ 
3200 & 1.25e-2 & 1.9425 & & 9.69e-6 & 2.9930 \\ 
6400 & 3.15e-3 & 1.9842 & & 1.21e-6 & 2.9965 \\  
\noalign{\smallskip}\hline
\end{tabular}
\label{tab:stiff_rk2}
\end{table}

\begin{table}[h!]
\renewcommand{\arraystretch}{1.1}
\scriptsize
\centering
\caption{Global error and order of accuracy for Example \ref{ex:stiff} by RK3 and RBF-RK3.}      
\begin{tabular}{clcrlcrlcrlc} 
\hline\noalign{\smallskip}
N & \multicolumn{2}{l}{RK3 I} & & \multicolumn{2}{l}{RBF-RK3 I} & & \multicolumn{2}{l}{RK3 II (a)} & & \multicolumn{2}{l}{RBF-RK3 II (a)} \\
    \cline{2-3}                   \cline{5-6}                       \cline{8-9}                        \cline{11-12}  
  & Error & Order             & & Error & Order                 & & Error & Order                  & & Error & Order \\
\noalign{\smallskip}\hline\noalign{\smallskip}
200  & 4.34e-2 & --     & & 2.94e-4  & --     & & 4.41e-2 & --     & & 3.14e-4  & --     \\  
400  & 5.85e-3 & 2.8907 & & 1.95e-5  & 3.9174 & & 5.91e-3 & 2.8973 & & 2.04e-5  & 3.9456 \\  
800  & 7.49e-4 & 2.9659 & & 1.25e-6  & 3.9587 & & 7.55e-4 & 2.9701 & & 1.30e-6  & 3.9729 \\
1600 & 9.46e-5 & 2.9856 & & 7.95e-8  & 3.9794 & & 9.51e-5 & 2.9879 & & 8.18e-8  & 3.9865 \\ 
3200 & 1.19e-5 & 2.9931 & & 5.00e-9  & 3.9913 & & 1.19e-5 & 2.9943 & & 5.12e-9  & 3.9984 \\ 
6400 & 1.49e-6 & 2.9966 & & 3.21e-10 & 3.9603 & & 1.50e-6 & 2.9972 & & 3.24e-10 & 3.9807 \\
\noalign{\smallskip}\hline\noalign{\smallskip}
N & \multicolumn{2}{l}{RK3 II (b)} & & \multicolumn{2}{l}{RBF-RK3 II (b)} & & \multicolumn{2}{l}{RK3 III (a)} & & \multicolumn{2}{l}{RBF-RK3 III (a)} \\
    \cline{2-3}                        \cline{5-6}                            \cline{8-9}                         \cline{11-12}  
  & Error & Order                  & & Error & Order                      & & Error & Order                   & & Error & Order \\
\noalign{\smallskip}\hline\noalign{\smallskip}
200  & 6.61e-2 & --     & & 3.93e-4 & --      & & 4.34e-2 & --     & & 3.75e-4  & --     \\  
400  & 9.10e-3 & 2.8596 & & 2.00e-5 &  4.2981 & & 5.83e-3 & 2.8985 & & 2.43e-5  & 3.9438 \\  
800  & 1.17e-3 & 2.9619 & & 1.30e-6 &  3.9429 & & 7.43e-3 & 2.9704 & & 1.55e-6  & 3.9720 \\
1600 & 1.47e-4 & 2.9854 & & 1.04e-4 & -6.3228 & & 9.37e-5 & 2.9880 & & 9.79e-8  & 3.9860 \\ 
3200 & 1.85e-5 & 2.9932 & & 3.35e-6 &  4.9563 & & 1.18e-5 & 2.9944 & & 6.14e-9  & 3.9963 \\ 
6400 & 2.32e-6 & 2.9967 & & 1.06e-7 &  4.9829 & & 1.47e-6 & 2.9972 & & 3.83e-10 & 4.0007 \\
\noalign{\smallskip}\hline\noalign{\smallskip}
N & \multicolumn{2}{l}{RK3 III (b)} & & \multicolumn{2}{l}{RBF-RK3 III (b)} & & \multicolumn{2}{l}{RK3 IV} & & \multicolumn{2}{l}{RBF-RK3 IV} \\
    \cline{2-3}                         \cline{5-6}                             \cline{8-9}                    \cline{11-12}  
  & Error & Order                   & & Error & Order                       & & Error & Order              & & Error & Order \\
\noalign{\smallskip}\hline\noalign{\smallskip}
200  & 6.75e-2 & --     & & 3.47e-4  & --     & & 4.77e-2 & --     & & 4.13e-4  & -- \\  
400  & 9.27e-3 & 2.8644 & & 4.46e-5  & 2.9629 & & 6.42e-3 & 2.8919 & & 2.70e-5  & 3.9375 \\  
800  & 1.19e-3 & 2.9653 & & 3.35e-6  & 3.7352 & & 8.20e-4 & 2.9692 & & 1.72e-6  & 3.9695 \\
1600 & 1.50e-4 & 2.9872 & & 2.24e-7  & 3.8988 & & 1.03e-4 & 2.9877 & & 1.09e-7  & 3.9855 \\ 
3200 & 1.88e-5 & 2.9942 & & 1.34e-8  & 4.0657 & & 1.30e-5 & 2.9942 & & 6.83e-9  & 3.9921 \\ 
6400 & 2.35e-6 & 2.9972 & & 2.07e-10 & 6.0193 & & 1.63e-6 & 2.9972 & & 4.23e-10 & 4.0125 \\
\noalign{\smallskip}\hline
\end{tabular}
\label{tab:stiff_rk3}
\end{table}

\begin{table}[h!]
\renewcommand{\arraystretch}{1.1}
\scriptsize
\centering
\caption{Global error and order of accuracy for Example \ref{ex:stiff} by RK4 and RBF-RK4.}      
\begin{tabular}{clcrlcrlc} 
\hline\noalign{\smallskip}
N & \multicolumn{2}{l}{RK4 I} & & \multicolumn{2}{l}{RBF RK4 I (+)} & & \multicolumn{2}{l}{RBF RK4 I (-)} \\
    \cline{2-3}                   \cline{5-6}                           \cline{8-9}
  & Error & Order             & & Error & Order                     & & Error & Order \\
\noalign{\smallskip}\hline\noalign{\smallskip}
200  & 6.19e-4  & --     & & 1.50e-4  & --      & & 6.66e-6  & --      \\  
400  & 3.98e-5  & 3.9582 & & 4.97e-6  &  4.9159 & & 6.17e-7  & 3.4335  \\   
800  & 2.52e-6  & 3.9801 & & 1.42e-7  &  5.1284 & & 4.52e-8  & 3.7698 \\ 
1600 & 1.59e-7  & 3.9902 & & 3.07e-9  &  5.5316 & & 3.04e-9  & 3.8924 \\  
3200 & 9.95e-9  & 3.9972 & & 3.11e-12 &  9.9452 & & 2.02e-10 & 3.9155 \\  
6400 & 6.20e-10 & 4.0037 & & 1.13e-11 & -1.8540 & & 1.66e-11 & 3.6036 \\   
\noalign{\smallskip}\hline\noalign{\smallskip}
N & \multicolumn{2}{l}{RK4 II} & & \multicolumn{2}{l}{RBF RK4 II (+)} & & \multicolumn{2}{l}{RBF RK4 II (-)} \\
    \cline{2-3}                    \cline{5-6}                            \cline{8-9}
  & Error & Order              & & Error & Order                      & & Error & Order \\
\noalign{\smallskip}\hline\noalign{\smallskip}
10  & 6.59e-4  & --     & & 6.98e-4  & --     & & 1.34e-6  & --     \\  
20  & 4.27e-5  & 3.9485 & & 2.57e-5  & 4.7650 & & 2.93e-7  & 2.1912 \\   
40  & 2.71e-6  & 3.9749 & & 8.62e-7  & 4.8969 & & 2.53e-8  & 3.5332 \\ 
80  & 1.71e-7  & 3.9877 & & 2.73e-8  & 4.9800 & & 1.80e-9  & 3.8099 \\  
160 & 1.07e-8  & 3.9936 & & 8.25e-10 & 5.0484 & & 1.31e-10 & 3.7802 \\  
320 & 6.70e-10 & 4.0018 & & 1.49e-11 & 5.7862 & & 1.35e-11 & 3.2863 \\   
\noalign{\smallskip}\hline
\end{tabular}
\label{tab:stiff_rk4}
\end{table}

\begin{example} \label{ex:nonseparable}
We conclude this section with an example, whose quadratic equation in the four-stage RBF Runge-Kutta method gives complex roots,
\begin{align*}
 \frac{du}{dt} &= \frac{2t^2-u}{t^2u-t}, \quad t \geqslant 1, \\
          u(1) &= 2,
\end{align*}
where the exact solution is $u(t) = \frac{1}{t} + \sqrt{\frac{1}{t^2}+4t-4}$ and the final time is $T=2$.

The global error at $t=T$ versus $N$ and the order of accuracy are shown in Tables \ref{tab:nonseparable_rk2}, \ref{tab:nonseparable_rk3} and \ref{tab:nonseparable_rk4}.
We see that the $s$-stage Runge-Kutta method maintains the expected order of accuracy.
The two-stage RBF Runge-Kutta method in Table \ref{tab:nonseparable_rk2} achieves third-order accuracy and exhibits better accuracy than the corresponding two-stage Runge-Kutta method for all $N$.
Some three-stage RBF Runge-Kutta methods attain fourth-order accuracy in Table \ref{tab:nonseparable_rk3}, where RBF-RK3 II (a) suffers from the order reduction while an order increase happens to RBF-RK3 II (b).
Most three-stage RBF Runge-Kutta methods are more accurate than the corresponding three-stage RBF Runge-Kutta methods at the same number of steps.
Any four-stage RBF Runge-Kutta method cannot reach order $5$ because the quadratic equations \eqref{eq:qd_eqn_I} and \eqref{eq:qd_eqn_II} produce complex roots at some step, leading to the complex approximation deviating far from the real exact solution and less accurate than the real approximation from the corresponding four-stage RBF Runge-Kutta method.
We also see that the $s$-stage RBF Runge-Kutta method performs better than most of the $(s\! +\! 1)$-stage Runge-Kutta methods in terms of accuracy, which is consistent with the observation in Example \ref{ex:stiff}.
\end{example}

\begin{table}[h!]
\renewcommand{\arraystretch}{1.1}
\scriptsize
\centering
\caption{Global error and order of accuracy for Example \ref{ex:nonseparable} by RK2 and RBF-RK2.}      
\begin{tabular}{clcrlc} 
\hline\noalign{\smallskip}
N & \multicolumn{2}{l}{RK2} & & \multicolumn{2}{l}{RBF RK2} \\ 
    \cline{2-3}                 \cline{5-6}
  & Error & Order           & & Error & Order \\
\noalign{\smallskip}\hline\noalign{\smallskip}
10  & 2.01e-3 & --     & & 1.16e-3 & --     \\  
20  & 4.42e-4 & 2.1885 & & 1.24e-4 & 3.2237 \\  
40  & 1.04e-4 & 2.0902 & & 1.43e-5 & 3.1182 \\
80  & 2.52e-5 & 2.0431 & & 1.71e-6 & 3.0597 \\ 
160 & 6.21e-6 & 2.0210 & & 2.09e-7 & 3.0299 \\ 
320 & 1.54e-6 & 2.0104 & & 2.59e-8 & 3.0150 \\  
\noalign{\smallskip}\hline
\end{tabular}
\label{tab:nonseparable_rk2}
\end{table}

\begin{table}[h!]
\renewcommand{\arraystretch}{1.1}
\scriptsize
\centering
\caption{Global error and order of accuracy for Example \ref{ex:nonseparable} by RK3 and RBF-RK3.}      
\begin{tabular}{clcrlcrlcrlc} 
\hline\noalign{\smallskip}
N & \multicolumn{2}{l}{RK3 I} & & \multicolumn{2}{l}{RBF-RK3 I} & & \multicolumn{2}{l}{RK3 II (a)} & & \multicolumn{2}{l}{RBF-RK3 II (a)} \\
    \cline{2-3}                   \cline{5-6}                       \cline{8-9}                        \cline{11-12}  
  & Error & Order             & & Error & Order                 & & Error & Order                  & & Error & Order \\
\noalign{\smallskip}\hline\noalign{\smallskip}
10  & 4.59e-5 & --     & & 2.99e-5  & --     & & 1.77e-4 & --     & & 1.81e-4 & --     \\  
20  & 1.91e-6 & 4.5877 & & 1.10e-6  & 4.7634 & & 2.25e-5 & 2.9728 & & 7.46e-6 & 4.6025 \\  
40  & 5.88e-7 & 1.6989 & & 4.85e-8  & 4.5041 & & 2.76e-6 & 3.0309 & & 7.56e-7 & 3.3028 \\
80  & 8.98e-8 & 2.7116 & & 2.50e-9  & 4.2769 & & 3.39e-7 & 3.0258 & & 2.26e-7 & 1.7399 \\ 
160 & 1.21e-8 & 2.8968 & & 1.42e-10 & 4.1421 & & 4.19e-8 & 3.0153 & & 1.13e-8 & 4.3252 \\ 
320 & 1.55e-9 & 2.9564 & & 8.44e-12 & 4.0713 & & 5.20e-9 & 3.0082 & & 5.00e-9 & 1.1770 \\
\noalign{\smallskip}\hline\noalign{\smallskip}
N & \multicolumn{2}{l}{RK3 II (b)} & & \multicolumn{2}{l}{RBF-RK3 II (b)} & & \multicolumn{2}{l}{RK3 III (a)} & & \multicolumn{2}{l}{RBF-RK3 III (a)} \\
    \cline{2-3}                        \cline{5-6}                            \cline{8-9}                         \cline{11-12}  
  & Error & Order                  & & Error & Order                      & & Error & Order                   & & Error & Order \\
\noalign{\smallskip}\hline\noalign{\smallskip}
10  & 2.33e-4 & --     & & 4.09e-2  & --     & & 2.53e-4 & --     & & 1.67e-5  & --     \\  
20  & 2.80e-5 & 3.0528 & & 8.01e-4  & 5.6737 & & 3.02e-5 & 3.0617 & & 8.33e-7  & 4.3258 \\  
40  & 3.38e-6 & 3.0541 & & 2.25e-5  & 5.1535 & & 3.62e-6 & 3.0611 & & 5.07e-8  & 4.0373 \\
80  & 4.12e-7 & 3.0342 & & 6.75e-7  & 5.0579 & & 4.41e-7 & 3.0371 & & 3.20e-9  & 3.9853 \\ 
160 & 5.08e-8 & 3.0189 & & 2.19e-8  & 4.9477 & & 5.44e-8 & 3.0200 & & 2.02e-10 & 3.9852 \\ 
320 & 6.31e-9 & 3.0099 & & 8.17e-10 & 4.7432 & & 6.75e-9 & 3.0104 & & 1.27e-11 & 3.9907 \\
\noalign{\smallskip}\hline\noalign{\smallskip}
N & \multicolumn{2}{l}{RK3 III (b)} & & \multicolumn{2}{l}{RBF-RK3 III (b)} & & \multicolumn{2}{l}{RK3 IV} & & \multicolumn{2}{l}{RBF-RK3 IV} \\
    \cline{2-3}                         \cline{5-6}                             \cline{8-9}                    \cline{11-12}  
  & Error & Order                   & & Error & Order                       & & Error & Order              & & Error & Order \\
\noalign{\smallskip}\hline\noalign{\smallskip}
10  & 3.20e-4 & --     & & 1.37e-3  & --     & & 3.30e-4 & --     & & 1.59e-4  & -- \\  
20  & 4.17e-5 & 2.9398 & & 2.20e-5  & 5.9530 & & 4.03e-5 & 3.0307 & & 7.96e-6  & 4.3199 \\  
40  & 5.27e-6 & 2.9829 & & 1.05e-7  & 7.7085 & & 4.94e-6 & 3.0300 & & 4.50e-7  & 4.1454 \\
80  & 6.61e-7 & 2.9948 & & 2.28e-8  & 2.2114 & & 6.10e-7 & 3.0184 & & 2.68e-8  & 4.0690 \\ 
160 & 8.27e-8 & 2.9982 & & 2.12e-9  & 3.4239 & & 7.57e-8 & 3.0100 & & 1.64e-9  & 4.0337 \\ 
320 & 1.03e-8 & 2.9993 & & 1.50e-10 & 3.8212 & & 9.43e-9 & 3.0052 & & 1.01e-10 & 4.0167 \\
\noalign{\smallskip}\hline
\end{tabular}
\label{tab:nonseparable_rk3}
\end{table}

\begin{table}[h!]
\renewcommand{\arraystretch}{1.1}
\scriptsize
\centering
\caption{Global error and order of accuracy for Example \ref{ex:nonseparable} by RK4 and RBF-RK4.}      
\begin{tabular}{clcrlcrlc} 
\hline\noalign{\smallskip}
N & \multicolumn{2}{l}{RK4 I} & & \multicolumn{2}{l}{RBF RK4 I (+)} & & \multicolumn{2}{l}{RBF RK4 I (-)} \\
    \cline{2-3}                   \cline{5-6}                           \cline{8-9}
  & Error & Order             & & Error & Order                     & & Error & Order \\
\noalign{\smallskip}\hline\noalign{\smallskip}
10  & 2.48e-6  & --     & & 4.87e-5  & --     & & 4.85e-5  & --     \\  
20  & 7.16e-8  & 5.1135 & & 2.19e-6  & 4.4766 & & 2.18e-6  & 4.4749 \\   
40  & 1.13e-8  & 2.6605 & & 1.16e-7  & 4.2402 & & 1.16e-7  & 4.2390 \\ 
80  & 9.08e-10 & 3.6408 & & 6.67e-9  & 4.1164 & & 6.67e-9  & 4.1162 \\  
160 & 6.27e-11 & 3.8557 & & 4.01e-10 & 4.0568 & & 4.01e-10 & 4.0567 \\  
320 & 4.10e-12 & 3.9348 & & 2.46e-11 & 4.0280 & & 2.46e-11 & 4.0280 \\   
\noalign{\smallskip}\hline\noalign{\smallskip}
N & \multicolumn{2}{l}{RK4 II} & & \multicolumn{2}{l}{RBF RK4 II (+)} & & \multicolumn{2}{l}{RBF RK4 II (-)} \\
    \cline{2-3}                    \cline{5-6}                            \cline{8-9}
  & Error & Order              & & Error & Order                      & & Error & Order \\
\noalign{\smallskip}\hline\noalign{\smallskip}
10  & 9.46e-7  & --     & & 6.02e-5  & --     & & 6.25e-5  & --     \\  
20  & 1.08e-7  & 3.1364 & & 1.94e-6  & 4.9517 & & 2.02e-6  & 4.9480 \\   
40  & 6.84e-9  & 3.9748 & & 8.68e-8  & 4.4855 & & 8.81e-8  & 4.5224 \\ 
80  & 4.21e-10 & 4.0227 & & 4.58e-9  & 4.2452 & & 4.59e-9  & 4.2609 \\  
160 & 2.60e-11 & 4.0147 & & 2.62e-10 & 4.1237 & & 2.63e-10 & 4.1280 \\  
320 & 1.62e-12 & 4.0088 & & 1.57e-11 & 4.0620 & & 1.57e-11 & 4.0631 \\   
\noalign{\smallskip}\hline
\end{tabular}
\label{tab:nonseparable_rk4}
\end{table}

\section{Concluding remarks} \label{sec:conclusion}
In this paper, we have constructed the RBF Runge-Kutta methods, depending on the Gaussian RBF Euler method, for the initial value problem.
By introducing the shape parameters to the $s$-stage Runge-Kutta methods with $s \leqslant 4$, it is possible to annihilate the $h^s$ term and obtain the order of $h^{s+1}$, so that the RBF Runge-Kutta methods achieve one order higher than the traditional Runge-Kutta ones analytically.
We have also proved the convergence of the proposed RBF Runge-Kutta methods.
The stability region of each RBF Runge-Kutta method is provided and compared with the one of its correspondent Runge-Kutta method.
Numerical experiments confirm the improved performance and great potential of RBF Runge-Kutta methods.

\section*{Acknowledgments}
This research is supported by National Research Foundation of Korea under the grant number 2021R1A2C3009648 and POSTECH Basic Science Research Institute under the NRF grant number NRF2021R1A6A1A1004294412.

\section*{Appendix} \label{appx}
\begin{align*}
 \A^3_1 &= \frac{1}{2} - b_2 c_2 - b_3 c_3, \\
 \A^3_2 &= \frac{1}{6} - \frac{1}{2} b_2 c^2_2 - \frac{1}{2} b_3 c^2_3, \\
 \B^3_2 &= \frac{1}{6} - a_{32} b_3 c_2, \\
 \A^3_3 &= \frac{1}{24} - \frac{1}{6} b_2 c^3_2 - \frac{1}{6} b_3 c^3_3, \\
 \B^3_3 &= \frac{1}{8} - a_{32} b_3 c_2 c_3, \\
 \D^3_3 &= \frac{1}{24} - \frac{1}{2} a_{32} b_3 c^2_2.
\end{align*}

\begin{align*}
 \A^4_1 &= \frac{1}{2} - b_2 c_2 - b_3 c_3 - b_4 c_4, \\ 
 \A^4_2 &= \frac{1}{6} - \frac{1}{2} b_2 c^2_2 - \frac{1}{2} b_3 c^2_3 - \frac{1}{2} b_4 c^2_4, \\
 \B^4_2 &= \frac{1}{6} - a_{32} b_3 c_2 - a_{42} b_4 c_2 - a_{43} b_4 c_3, \\
 \A^4_3 &= \frac{1}{24} - \frac{1}{6} b_2 c^3_2 - \frac{1}{6} b_3 c^3_3 - \frac{1}{6} b_4 c^3_4, \\
 \B^4_3 &= \frac{1}{8} - a_{32} b_3 c_2 c_3 - a_{42} b_4 c_2 c_4 - a_{43} b_4 c_3 c_4, \\ 
 \D^4_3 &= \frac{1}{24} - \frac{1}{2} a_{32} b_3 c^2_2 - \frac{1}{2} a_{42} b_4 c^2_2 - \frac{1}{2} a_{43} b_4 c^2_3, \\
 \E^4_3 &= \frac{1}{24} - a_{32} a_{43} b_4 c_2, \\
 \A^4_4 &= \frac{1}{120} - \frac{1}{24} b_2 c_2^4 - \frac{1}{24} b_3 c_3^4 - \frac{1}{24} b_4 c_4^4, \\
 \B^4_4 &= \frac{1}{20} - \frac{1}{2} a_{32} b_3 c_2 c_3^2 - \frac{1}{2} a_{42} b_4 c_2 c_4^2 - \frac{1}{2} a_{43} b_4 c_3 c_4^2, \\
 \D^4_4 &= \frac{1}{30} - \frac{1}{2} a_{32} b_3 c_2^2 c_3 - \frac{1}{2} a_{42} b_4 c_2^2 c_4 - \frac{1}{2} a_{43} b_4 c_3^2 c_4, \\
 \E^4_4 &= \frac{1}{120} - \frac{1}{6} a_{32} b_3 c_2^3 - \frac{1}{6} a_{42} b_4 c_2^3 - \frac{1}{6} a_{43} b_4 c_3^3, \\
 \F^4_4 &= \frac{7}{120} - a_{32} a_{43} b_4 c_2 c_3 - a_{32} a_{43} b_4 c_2 c_4, \\
 \G^4_4 &= \frac{1}{120} - \frac{1}{2} a_{32} a_{43} b_4 c_2^2, \\
 \HH^4_4 &= \frac{1}{40} - \frac{1}{2} a_{32}^2 b_3 c_2^2 - \frac{1}{2} a_{42}^2 b_4 c_2^2 - \frac{1}{2} a_{43}^2 b_4 c_3^2 - a_{42} a_{43} b_4 c_2 c_3.
\end{align*}


\end{document}